\documentclass{amsart}

\usepackage{amsmath,amsthm,amssymb,amsfonts}
\usepackage{hyperref}

\input xy
\xyoption{all}
\usepackage{enumerate}
\usepackage{pgf,tikz}
\usetikzlibrary{matrix}
\usepackage{mathrsfs}
\usetikzlibrary{arrows}

\theoremstyle{plain}
\newtheorem{thm}{Theorem}[section]

\newtheorem{cor}[thm]{Corollary}

\theoremstyle{definition}

\newtheorem{rmk}[thm]{Remark}
\newtheorem{exm}[thm]{Example}
\newtheorem{exe}[thm]{Exercise}
\newtheorem{sol}[thm]{Solution}

\newcommand{\Hom}{\operatorname{Hom}}

\newcommand{\Ext}{\operatorname{Ext}}

\newcommand{\Img}{\operatorname{Im}}

\newcommand{\total}{\operatorname{Total}}

\makeatletter
\@tempcnta=\@ne
\def\@nameedef#1{\expandafter\edef\csname #1\endcsname}
\loop\ifnum\@tempcnta<27
  \@nameedef{C\@Alph\@tempcnta}{\noexpand\mathcal{\@Alph\@tempcnta}}
  \advance\@tempcnta\@ne
\repeat

\makeatletter
\@tempcnta=\@ne
\def\@nameedef#1{\expandafter\edef\csname #1\endcsname}
\loop\ifnum\@tempcnta<27
  \@nameedef{B\@Alph\@tempcnta}{\noexpand\mathbb{\@Alph\@tempcnta}}
  \advance\@tempcnta\@ne
\repeat

\makeatletter
\@tempcnta=\@ne
\def\@nameedef#1{\expandafter\edef\csname #1\endcsname}
\loop\ifnum\@tempcnta<27
  \@nameedef{F\@Alph\@tempcnta}{\noexpand\mathsf{\@Alph\@tempcnta}}
  \advance\@tempcnta\@ne
\repeat

\begin{document}
\title{Spectral sequences via examples}

\author{Antonio D\'iaz Ramos}
\email{adiazramos@uma.es}
\address{Departamento de \'Algebra, Geometr\'ia y Topolog\'ia, 
	Universidad de M\'alaga, 
	Apdo correos 59, 
	29080 M\'alaga, 
	Spain}
\thanks{\rm Author supported by MICINN grant RYC-2010-05663 and partially supported by FEDER-MCI grant MTM2013-41768-P and Junta de Andaluc{\'\i}a grant FQM-213.}

\date{\today}
\begin{abstract}
These are lecture notes for a short course about spectral sequences that was held at M\'alaga, October 18--20 (2016), during the ``Fifth Young Spanish Topologists Meeting''. The approach is to illustrate the basic notions via fully computed examples arising from Algebraic Topology and Group Theory. 
\end{abstract}
\maketitle

\begin{center}
\includegraphics[scale=0.15]{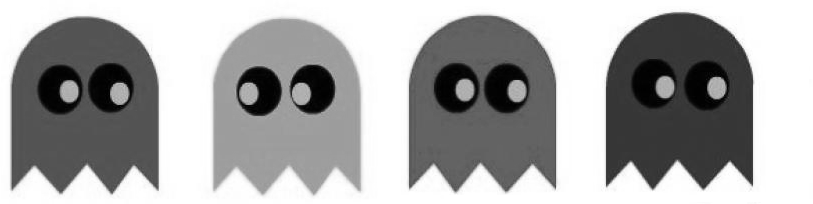}
\end{center}

\section{Foreword}

These notes avoid most of the technicalities about spectral sequences and some of the definitions are  simplified. The aim is to introduce the basic notions involved when working with spectral sequences with the minimum necessary background. We go along the way bringing in notions and illustrating them by concrete computations. The target reader is someone not acquainted with this topic. 
For full details the reader should consult \cite{McCleary}, \cite{Weibel} or \cite{MacLane}.

The following list describes the contents of the different sections of these notes:
\begin{enumerate}
\item[\eqref{ss2:ssandfil}] \textbf{Spectral sequences and filtrations.} We explain the intimate relation between spectral sequences, filtrations, quotients and (co)homology.
\item[\eqref{ss3:diffandconv}] \textbf{Differentials and convergence.} We give definitions of what a (homological type) spectral sequence is and what is convergence in this setup. We also talk about extension problems. 
\item[\eqref{ss4:cohotype}] \textbf{Cohomological type spectral sequences.} We introduce cohomological type spectral sequences and say a few more words about the extension problem.
\item[\eqref{ss5:algebra}] \textbf{Additional structure: algebra.}  We equip everything with products, yielding graded algebras and bigraded algebras. We start discussing the lifting problem.
\item[\eqref{ss6:lifting}] \textbf{Lifting problem.} We go in detail about the lifting problem and discuss several particular cases.
\item[\eqref{ss7:edgemorphisms}] \textbf{Edge morphisms.} After introducing the relevant notions, a couple of examples where the edge morphisms play a central role are fully described.
\item[\eqref{ss8:diffsssametarget}] \textbf{Different spectral sequences with same target.} This section is included to make the reader aware of this fact.
\item[\eqref{ss9:glimpse}] \textbf{A glimpse into the black box.}  We sketch the technicalities that have been avoided in the rest of the sections. In particular, we provide some detail on the construction of a spectral sequence from a filtration of a chain complex.
\item[\eqref{section:AppendixI}] \textbf{Appendix I.} Here, precise statements for the following spectral sequences are given: Serre spectral sequence of a fibration, Lyndon-Hochschild-Serre spectral sequence of a short exact sequence of groups, Atiyah-Hirzebruch spectral sequence of a fibration. All examples in these notes are based on these spectral sequences. Also, sketch of proofs of the two first spectral sequences are included and the filtrations employed are described.
\item[\eqref{section:AppendixII}] \textbf{Appendix II.} We list some important results in Topology whose proofs involve spectral sequences and we briefly comment on the role they play in the proofs.
\item[\eqref{ss:Solutions}] \textbf{Solutions.} Answers to the proposed exercises. Some of them are just a reference to where to find a detailed solution.
\end{enumerate}
We denote by $R$ the base (commutative with unit) ring over which we work. 

\vspace{2cm}

\textbf{Acknowledgements:} I am grateful to Oihana Garaialde, David M\'endez and Jose Manuel Moreno for reading a preliminary version of these notes. Also, thanks to Urtzi Buijs for its nice version of the bridge-river-tree picture using the TikZ package.

\newpage
\section{Spectral sequences and filtrations}
\label{ss2:ssandfil}
A spectral sequence is a tool to compute the homology of a chain complex. They are employed whenever we are unable to directly calculate this homology but instead we have a filtration of the chain complex. Each such filtration gives rise to a spectral sequence that may or may not help in determining the homology we are interested in.

A spectral sequence consists of a sequence of pages, $E^2,E^3,\ldots, E^\infty$, equipped with differentials, $d_2,d_3,\ldots$, in such a way that $E^{n+1}$ is the homology of $E^n$. If $H$ is the homology we wish to find out, we expect to gain from $E^2$ to $E^\infty$ using the differentials and then to determine $H$ from $E^\infty$.

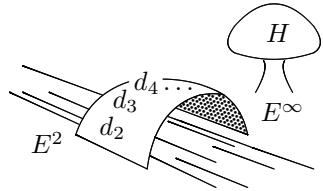
\begin{figure}[h]
\label{figure:bridgeandtree}
\begin{picture}(100,65)(0,-20)
\qbezier(0,0)(5,15)(20,20)
%\qbezier(20,20)(27,22)(33,21)
\qbezier(33,21)(40,21)(47,19)
\qbezier(47,15)(30,10)(27,-13)
\qbezier(47,15)(53,16)(56,15)
\qbezier(47,19)(53,17)(56,15)
\qbezier(56,15)(63,10)(65,0)

\qbezier(0,0)(0,0)(27,-13)
\qbezier(38,10)(65,0)(65,0)

%RIO
\qbezier(28,-8)(28,-8)(53,-20)
%\qbezier(33,-8)(33,-8)(60,-21)
%\qbezier(40,-6)(40,-6)(67,-17)
\qbezier(38,-9)(38,-9)(65,-20)
\qbezier(40,5)(40,5)(67,-5)
\qbezier(35,6.8)(40,5)(40,5)
\qbezier(67,-5)(67,-5)(90,-13)
%\qbezier(57,-5)(57,-5)(80,-13)
\qbezier(52,-4.5)(52,-4.5)(75,-13)

\qbezier(30,-1)(30,-1)(53,-10.5)
\qbezier(35,5)(35,5)(60,-4.5)

\qbezier(-25,16)(-25,16)(1,4)
\qbezier(-20,26)(-20,26)(10,15)

\qbezier(-20,16)(-20,16)(2.7,6)
\qbezier(-10,16)(-10,16)(5,9.6)

\put(-17,-7){$E^2$}
\put(9,0){$d_2$}
\put(14,10){$d_3$}
\put(22,17){$d_4$}
\put(33,17){$\dots $}

\put(70,5){$E^\infty$}

\qbezier(68,16)(73,25)(71,28)
\qbezier(82,16)(77,25)(79,28)

%COPA
\qbezier(60,31)(75,26)(90,31)
\qbezier(60,31)(55, 34)(60,41)
\qbezier(90,31)(95,34)(90, 41)
\qbezier(60, 41)(75,58)(90, 41)

%SOMBRA
\multiput(38.5,8)(1.5,-0.55){17}{$\cdot$}
\multiput(39.5,9)(1.5,-0.55){16}{$\cdot$}
\multiput(40.5,10)(1.5,-0.55){15}{$\cdot$}
\multiput(43,10.45)(1.5,-0.55){13}{$\cdot$}
\multiput(44,11.45)(1.5,-0.55){12}{$\cdot$}
\multiput(46.5,11.90)(1.5,-0.55){9}{$\cdot$}
\multiput(49,12.35)(1.5,-0.55){7}{$\cdot$}
\multiput(53,12.25)(1.5,-0.55){2}{$\cdot$}
\put(72,35){$H$}
\end{picture}
\caption{You want to cross the bridge and then to climb the tree. Don't fall into the river of madness!}
\end{figure}

%\begin{figure}[h]
%\label{figure:bridgeandtree}
%\begin{tikzpicture}[scale=0.2,line cap=round,line join=round,>=triangle 45,x=1.0cm,y=1.0cm]
%\clip(-2.,3.) rectangle (18.,15.);
%\draw [line width=2.pt] (2.,4.)-- (2.,9.);
%\draw [line width=2.pt] (2.,4.)-- (4.,4.);
%\draw [line width=2.pt] (10.,4.)-- (12.,4.);
%\draw [line width=2.pt] (12.,4.)-- (12.,9.);
%\draw [line width=2.pt] (12.,9.)-- (2.,9.);
%\draw [shift={(7.,4.)},line width=2.pt]  plot[domain=0.:3.141592653589793,variable=\t]({1.*3.*cos(\t r)+0.*3.*sin(\t r)},{0.*3.*cos(\t r)+1.*3.*sin(\t r)});
%\draw (-2.24,7.18) node[anchor=north west] {$E^2$};
%\draw (13.46,7.26) node[anchor=north west] {$E^\infty$};
%\draw [line width=2.pt] (14.02,7.48)-- (14.02,12.46);
%\draw [line width=2.pt] (16.02,7.48)-- (16.02,12.48);
%\draw [line width=2.pt] (14.02,12.)-- (16.02,12.);
%\draw [line width=2.pt] (14.02,11.)-- (16.02,11.);
%\draw [line width=2.pt] (14.02,10.)-- (16.02,10.);
%\draw [line width=2.pt] (14.02,9.)-- (16.02,9.);
%\draw [line width=2.pt] (14.02,8.)-- (16.02,8.);
%\draw (13.8,15.38) node[anchor=north west] {$H$};
%\draw (2.,7.6) node[anchor=north west] {$d_2$};
%\draw (3.3,8.9) node[anchor=north west] {$d_3$};
%\draw (5.16,9.3) node[anchor=north west] {$d_4$};
%\draw (7.08,9.2) node[anchor=north west] {$d_5$};
%\draw (8.72,7.5) node[anchor=north west] {$...$};
%\end{tikzpicture}
%\caption{You want to cross the bridge and then to climb the tree. Don't fall into the river of madness!}
%\end{figure}

A typical theorem about spectral sequences tells you only about $H$ and $E^2$. The filtration being utilised to construct the spectral sequence usually reveals a deep fact about the context, and the description of the differentials is habitually missing, as they are generally unknown. In spite of these apparent limitations, spectral sequences are very useful, see Appendix II \ref{section:AppendixII} for some deep applications to Topology.

You might wonder whatever happened to $E^0$ and $E^1$. In fact, spectral sequences start at $E^0$. As these two first pages may depend on a choice or may not have a nice description,  they are often, but not always, ignored in the statement of a theorem about a spectral sequence. 

\begin{exm}\label{example:semisimplicialmodelofS^2}
Consider a (semi-simplicial) model of the $2$-sphere $S^2$ with vertices $\{a,b,c\}$, edges $\{A,B,C\}$ and solid triangles $\{P,Q\}$ and with inclusions as follows:
\begin{minipage}{0.45\textwidth}
\[
\xymatrix@=20pt{
P\ar@{-}[]+R;[rrdd]+L\ar@{-}[]+R;[rr]+L\ar@{-}[]+R;[rrd]+L&& A \ar@{-}[]+R;[rr]+L\ar@{-}[]+R;[rrd]+L&& a\\
	&& B\ar@{-}[]+R;[rru]+L\ar@{-}[]+R;[rrd]+L&&b\\
Q\ar@{-}[]+R;[rruu]+L\ar@{-}[]+R;[rru]+L\ar@{-}[]+R;[rr]+L&& C \ar@{-}[]+R;[rr]+L\ar@{-}[]+R;[rru]+L&& c
}
\]
\end{minipage}
\begin{minipage}{0.25\textwidth}
\definecolor{zzttqq}{rgb}{0.26666666666666666,0.26666666666666666,0.26666666666666666}
\definecolor{qqqqff}{rgb}{0.3333333333333333,0.3333333333333333,0.3333333333333333}
\begin{tikzpicture}[scale=0.4,line cap=round,line join=round,>=triangle 45,x=1.0cm,y=1.0cm]
\clip(7.,1.) rectangle (14.,14.);
\fill[color=zzttqq,fill=zzttqq,fill opacity=0.1] (12.445958331289594,12.82878580241334) -- (8.043187170820223,11.379955178578085) -- (13.262220865521774,9.712147035345012) -- cycle;
\fill[color=zzttqq,fill=zzttqq,fill opacity=0.1] (12.328565475539177,4.467382045067803) -- (8.,3.) -- (13.219033694701556,1.3321918567669275) -- cycle;
\draw (8.,7.)-- (13.219033694701551,5.3321918567669275);
\draw (13.219033694701551,5.3321918567669275)-- (12.328565475539172,8.467382045067803);
\draw (12.328565475539172,8.467382045067803)-- (8.,7.);
\draw (8.043187170820223,11.379955178578085)-- (13.262220865521774,9.712147035345012);
\draw (13.262220865521774,9.712147035345012)-- (12.445958331289594,12.82878580241334);
\draw (12.445958331289594,12.82878580241334)-- (8.043187170820223,11.379955178578085);
\draw (8.,3.)-- (13.219033694701556,1.3321918567669275);
\draw (13.219033694701556,1.3321918567669275)-- (12.328565475539177,4.467382045067803);
\draw (12.328565475539177,4.467382045067803)-- (8.,3.);
\draw [->,dash pattern=on 2pt off 2pt] (8.043187170820225,11.379955178578083) -- (8.024635749587674,7.428502456045029);
\draw [->,dash pattern=on 2pt off 2pt] (13.262220865521773,9.712147035345012) -- (13.367445064561949,6.14845439099911);
\draw [->,dash pattern=on 2pt off 2pt] (12.445958331289594,12.828785802413336) -- (12.365668318004271,9.09813036697449);
\draw [->,dash pattern=on 2pt off 2pt] (12.328565475539179,4.467382045067803) -- (12.291462633074072,7.4841567197426775);
\draw [->,dash pattern=on 2pt off 2pt] (8.,3.) -- (8.006084328355124,6.259762918394406);
\draw [->,dash pattern=on 2pt off 2pt] (13.219033694701556,1.3321918567669273) -- (13.2375851159341,4.664340692395144);
\begin{scriptsize}
\draw [fill=qqqqff] (8.,7.) circle (2.5pt);
\draw (8.302907068075918,7.3357453498822816) node {$a$};
\draw [fill=qqqqff] (13.219033694701551,5.3321918567669275) circle (2.5pt);
\draw (13.460202170724695,5.666117438952821) node {$b$};
\draw[color=black] (10.566180458446963,5.8) node {$A$};
\draw [fill=qqqqff] (12.328565475539172,8.467382045067803) circle (2.5pt);
\draw (12.569733951562316,8.801307627253696) node {$c$};
\draw[color=black] (13,7.5) node {$C$};
\draw[color=black] (10.120946348865774,8.17055930534701) node {$B$};
\draw (11.289685886516397,11.472712284740833) node {$P$};
\draw (11.234031622818748,3.087469887628431) node {$Q$};
\end{scriptsize}
\end{tikzpicture}
\end{minipage}
\begin{minipage}{0.25\textwidth}
\definecolor{xdxdff}{rgb}{0.3333333333333333,0.3333333333333333,0.3333333333333333}
\begin{tikzpicture}[scale=0.25,line cap=round,line join=round,>=triangle 45,x=1.0cm,y=1.0cm]
\clip(3.,2.) rectangle (16.,13.);
\draw [fill=black,fill opacity=0.1] (9.74,7.56) circle (5.2288813335167585cm);
\draw [shift={(9.66,11.44)}] plot[domain=3.7824208010115905:5.678013499634892,variable=\t]({1.*6.423057762301059*cos(\t r)+0.*6.423057762301059*sin(\t r)},{0.*6.423057762301059*cos(\t r)+1.*6.423057762301059*sin(\t r)});
\draw (9.92,10.3) node[anchor=north west] {P};
\draw (9,4.5) node[anchor=north west] {Q};
\draw (9.76,6.8) node[anchor=north west] {A};
\begin{scriptsize}
\draw [fill=xdxdff] (5.354356892643101,6.673753211390967) circle (2.5pt);
\draw(5.68,7.04) node {$a$};
\draw [fill=xdxdff] (13.82856110438738,6.553409231698999) circle (2.5pt);
\draw (13.8,7.2) node {$b$};
\end{scriptsize}
\end{tikzpicture}
\end{minipage}

The homology $H$ of the following chain complex $(C,d)$ gives the integral homology of the sphere $S^2$:
\[
\xymatrix@C=10pt@R=0pt{
0\ar[r] &\BZ\{P,Q\}\ar[r]^d& \BZ\{A,B,C\}\ar[r]^d & \BZ\{a,b,c\}\ar[r]& 0.\\
& d(P)=C-B+A 		 & d(A)=b-a\\
& d(Q)=C-B+A 		 & d(B)=c-a\\
			&	     & d(C)=c-b
}
\]
Instead of directly calculating $H$, we instead consider the following filtration:
\[
\xymatrix@C=10pt@R=0pt{
0\ar[r] &\BZ\{P,Q\}\ar[r]& \BZ\{A,B,C\}\ar[r] & \BZ\{a,b,c\}\ar[r]& 0\\
0\ar[r] &0\ar[r]& \BZ\{A,B\}\ar[r] & \BZ\{a,b,c\}\ar[r]& 0\\
0\ar[r] &0\ar[r]& \BZ\{A\}\ar[r] & \BZ\{a,b\}\ar[r]& 0.\\
}
\]
As you can see, the differential restricts to each row and hence makes each row into a chain complex. Moreover, each row is contained in the row above and hence we may consider the quotient of successive rows. We call this quotient $E^0$ and it inherits a differential $d_0$ from $d$:
\[
\xymatrix@C=10pt@R=0pt{
0\ar[r] &\BZ\{P,Q\}\ar[r]& \BZ\{C\}\ar[r] & 0\ar[r]& 0&&d_0(P)=C,d_0(Q)=C\\
0\ar[r] &0\ar[r]& \BZ\{B\}\ar[r] & \BZ\{c\}\ar[r]& 0&&d_0(B)=c\\
0\ar[r] &0\ar[r]& \BZ\{A\}\ar[r] & \BZ\{a,b\}\ar[r]& 0&&d_0(A)=b-a.\\
}
\]
If we compute the homology with respect to $d_0$ we obtain $E^1$:
\[
\xymatrix@C=10pt@R=0pt{
0\ar[r] &\BZ\{P-Q\}\ar[r]& 0 \ar[r] & 0\ar[r]& 0\\
0\ar[r] &0\ar[r]& 0\ar[r] & 0\ar[r]& 0\\
0\ar[r] &0\ar[r]& 0\ar[r] & \BZ\{\bar a\}\ar[r]& 0.\\
}
\]
Clearly, the differential $d_1$, which is also induced by $d$, is zero. This implies that $E^2=E^1$. For the same reason, $d_n=0$ for all $n\geq 2$ and $E^\infty=\ldots=E^3=E^2=E^1$. The two classes $P-Q$ and $\bar a$ correspond respectively to generators of the homology groups $H_2(S^2;\BZ)=\BZ$ and $H_0(S^2;\BZ)=\BZ$.
\end{exm}

\begin{exe}
Cook up a different filtration to that in Example \ref{example:semisimplicialmodelofS^2} and compute $E^0$ and $E^1$. Obtain again the integral homology of $S^2$.
\end{exe}

We can learn much from this example. First fact is that, in general,
\[
\text{$E^0$ is the quotient of $C$ by the filtration,}
\]
and $d_0$ is induced by $d$. Then $E^1$ is the homology of $(E^0,d_0$) and $d_1$ is induced by $d$ again. In fact, this is the general rule to pass from $E^n$ to $E^{n+1}$:
\[
\text{$E^{n+1}$ is the homology of $(E^n,d_n)$ and $d_{n+1}$ is induced from $d$.}
\]
Although there exists a closed expression for $d_{n+1}$ in terms of $d$ and $C$, it is so complicated that, in practice, you are only able to fully describe $(E^0,d_0)$, $(E^1,d_1)$ and $E^2$. Then $E^2$ is what you write in the statement of your theorem about your spectral sequence. You are seldom capable of figuring out what $d_2, d_3,\ldots$ are. 

This gives you some hints about the bridge in Figure \ref{figure:bridgeandtree}. The next example provides information on how the tree in that figure looks like.

\begin{exm}\label{example:semisimplicialmodelofS^2Einftyexplained}
Recall the filtration of chain complexes we gave in Example \ref{example:semisimplicialmodelofS^2}:
\[
\xymatrix@C=10pt@R=0pt{
0\ar[r] &\BZ\{P,Q\}\ar[r]& \BZ\{A,B,C\}\ar[r] & \BZ\{a,b,c\}\ar[r]& 0\\
0\ar[r] &0\ar[r]& \BZ\{A,B\}\ar[r] & \BZ\{a,b,c\}\ar[r]& 0\\
0\ar[r] &0\ar[r]& \BZ\{A\}\ar[r] & \BZ\{a,b\}\ar[r]& 0.\\
}
\]
Now, instead of taking quotient followed by homology, we do it the other way around, i.e.,  we first take homology,
\[
\xymatrix@C=10pt@R=0pt{
0\ar[r] &\BZ\{P-Q\}\ar[r]& 0\ar[r] & \BZ\{\bar a\}\ar[r]& 0\\
0\ar[r] &0\ar[r]& 0\ar[r] & \BZ\{\bar a\}\ar[r]& 0\\
0\ar[r] &0\ar[r]& 0\ar[r] & \BZ\{\bar a\}\ar[r]& 0,\\
}
\]
and then we take quotient by successive rows,
\[
\xymatrix@C=10pt@R=0pt{
0\ar[r] &\BZ\{P-Q\}\ar[r]& 0 \ar[r] & 0\ar[r]& 0\\
0\ar[r] &0\ar[r]& 0\ar[r] & 0\ar[r]& 0\\
0\ar[r] &0\ar[r]& 0\ar[r] & \BZ\{\bar a\}\ar[r]& 0.\\
}
\]
As you can observe, we have obtained again $E^\infty=E^1$ of Example \ref{example:semisimplicialmodelofS^2}.
\end{exm}

What happened in the above example is not a coincidence and, in general,
\[
\text{$E^\infty$ is the quotient of $H$ by the induced filtration.}
\]
We finish this section with an enriched version of our initial picture:

\begin{figure}[h]
\begin{picture}(100,95)(20,-40)
\put(-105,50){something you probably don't know}
\put(20,35){$\Big\downarrow$}
\qbezier(0,0)(5,15)(20,20)
%\qbezier(20,20)(27,22)(33,21)
\qbezier(33,21)(40,21)(47,19)
\qbezier(47,15)(30,10)(27,-13)
\qbezier(47,15)(53,16)(56,15)
\qbezier(47,19)(53,17)(56,15)
\qbezier(56,15)(63,10)(65,0)

\qbezier(0,0)(0,0)(27,-13)
\qbezier(38,10)(65,0)(65,0)

%RIO
\qbezier(28,-8)(28,-8)(53,-20)
%\qbezier(33,-8)(33,-8)(60,-21)
%\qbezier(40,-6)(40,-6)(67,-17)
\qbezier(38,-9)(38,-9)(65,-20)
\qbezier(40,5)(40,5)(67,-5)
\qbezier(35,6.8)(40,5)(40,5)
\qbezier(67,-5)(67,-5)(90,-13)
%\qbezier(57,-5)(57,-5)(80,-13)
\qbezier(52,-4.5)(52,-4.5)(75,-13)

\qbezier(30,-1)(30,-1)(53,-10.5)
\qbezier(35,5)(35,5)(60,-4.5)

\qbezier(-25,16)(-25,16)(1,4)
\qbezier(-20,26)(-20,26)(10,15)

\qbezier(-20,16)(-20,16)(2.7,6)
\qbezier(-10,16)(-10,16)(5,9.6)

\put(-85,-7){$E^0\stackrel{d_0}\longrightarrow E^1\stackrel{d_1}\longrightarrow$}
\put(-82,-22){$\Big\uparrow$}
\put(-100,-40){quotient of $C$ by filtration}
\put(-17,-7){$E^2$}
\put(-14,-17){$\uparrow$}
\put(-30,-28){something you know}
\put(9,0){$d_2$}
\put(14,10){$d_3$}
\put(22,17){$d_4$}
\put(33,17){$\dots $}

\put(70,5){$E^\infty\leftarrow$quotient of $H$ by induced filtration}

\qbezier(68,16)(73,25)(71,28)
\qbezier(82,16)(77,25)(79,28)

%COPA
\qbezier(60,31)(75,26)(90,31)
\qbezier(60,31)(55, 34)(60,41)
\qbezier(90,31)(95,34)(90, 41)
\qbezier(60, 41)(75,58)(90, 41)

%SOMBRA
\multiput(38.5,8)(1.5,-0.55){17}{$\cdot$}
\multiput(39.5,9)(1.5,-0.55){16}{$\cdot$}
\multiput(40.5,10)(1.5,-0.55){15}{$\cdot$}
\multiput(43,10.45)(1.5,-0.55){13}{$\cdot$}
\multiput(44,11.45)(1.5,-0.55){12}{$\cdot$}
\multiput(46.5,11.90)(1.5,-0.55){9}{$\cdot$}
\multiput(49,12.35)(1.5,-0.55){7}{$\cdot$}
\multiput(53,12.25)(1.5,-0.55){2}{$\cdot$}
\put(72,35){$H\longleftarrow$something you want to know}
\end{picture}
\caption{A spectral sequence measures the lack of commutativity of the operations taking homology and taking quotient.}
\end{figure}
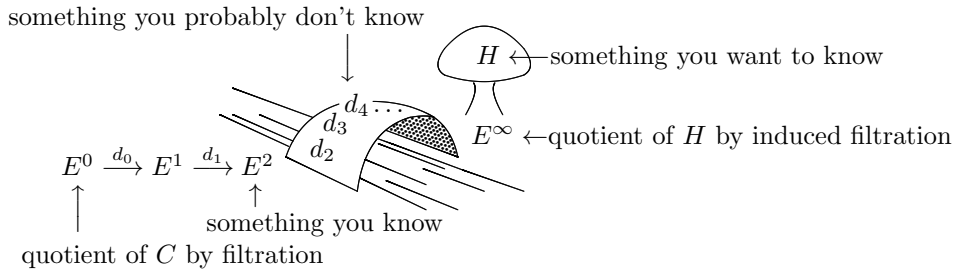

\section{Differentials and convergence}
\label{ss3:diffandconv}
There are homological type and cohomological type spectral sequences. The former type  consists of a sequence of pages and differentials $\{E^r,d_r\}_{r\geq 2}$, such that $E^{r+1}$ is the homology of $(E^r,d_r)$. Each page is a differential bigraded module over the base ring $R$. This means that at each position $(p,q)$, there is an $R$-module $E^r_{p,q}$. Moreover, the differential $d_r$ squares to zero and has bidegree $(-r,r-1)$. This last bit means that it maps $d_r\colon E^r_{p,q}\to E^r_{p-r,q+r-1}$, as in these pictures:

\begin{minipage}{0.3\textwidth}
\begin{align*}
&\xymatrix@=0pt{
 \bullet & \bullet & \bullet & \bullet & \bullet \\
 \bullet & \bullet & \bullet & \bullet & \bullet \\
 \bullet & \bullet & \ar[]+C;[llu]+C\bullet & \bullet & \bullet \\
 \bullet & \bullet & \bullet & \bullet & \bullet \\
 \bullet & \bullet & \bullet & \bullet & \bullet \\
}\\
&\hspace{10pt}(E^2,d_2)
\end{align*}
\end{minipage}
\begin{minipage}{0.3\textwidth}
\begin{align*}
&\xymatrix@=0pt{
 \bullet & \bullet & \bullet & \bullet & \bullet \\
 \bullet & \bullet & \bullet & \bullet & \bullet \\
 \bullet & \bullet & \bullet & \bullet & \bullet \\
 \bullet & \bullet & \bullet & \ar[]+C;[llluu]+C\bullet & \bullet \\
 \bullet & \bullet & \bullet & \bullet & \bullet \\
}\\
&\hspace{10pt}(E^3,d_3)
\end{align*}
\end{minipage}
\begin{minipage}{0.3\textwidth}
\begin{align*}
&\xymatrix@=0pt{
 \bullet & \bullet & \bullet & \bullet & \bullet \\
 \bullet & \bullet & \bullet & \bullet & \bullet \\
 \bullet & \bullet & \bullet & \bullet & \bullet \\
 \bullet & \bullet & \bullet & \bullet &  \ar[]+C;[lllluuu]+C\bullet \\
 \bullet & \bullet & \bullet & \bullet & \bullet \\
}\\
&\hspace{10pt}(E^4,d_4)
\end{align*}
\end{minipage}

One does not need to know neither $(E^0,d_0)$ nor $(E^1,d_1)$ to start computing with a spectral sequence. You can proceed even if you do not know $C$ or the filtration that gave rise to the spectral sequence. If you know $H$ and you are told that the spectral sequence \emph{converges to $H$}, written as $E^2_{p,q}\Rightarrow H_{p+q}$, this roughly means that:
\begin{enumerate}[(a)]
\item \label{enumerate:convergencestabilizes} the value at each position $(p,q)$ stabilizes after a finite number of steps, and the stable value is denoted by $E^\infty_{p,q}$, and that
\item \label{enumerate:convergencefiniteextensions} you can recover $H$ from $E^\infty$.
\end{enumerate} 

Before unfolding details about point \eqref{enumerate:convergencefiniteextensions}, we work out an example.

\begin{exm}\label{example:HopffibrationS1S3S2}
The sphere of dimension $S^n$ is a Moore space $M(n,\BZ)$, i.e., its homology is given by $H_i(S^n;\BZ)=\BZ$ if $i=n,0$ and $0$ otherwise. Moreover, the spheres of dimensions $1$, $3$ and $2$ nicely fit into the Hopf fibration:

\begin{minipage}{0.4\textwidth}
\[
\xymatrix{
S^1\ar[r]& S^3\ar[r]& S^2}
\]
\end{minipage}
\begin{minipage}{0.4\textwidth}
\begin{tikzpicture}[scale=0.4,line cap=round,line join=round,>=triangle 45,x=1.0cm,y=1.0cm]
\clip(2,3) rectangle (17,8);
\draw [fill=black,fill opacity=0.25] (15,5.479480519480518) circle (1.6775101241535375cm);
\draw(5.144069264069266,5.566060606060605) circle (1.6699699429038262cm);
\draw (6.9,5.929696969696968) node[anchor=north west] {$\longrightarrow \mathbb{R}^3\sqcup *\longrightarrow $};
\draw [shift={(15.1,9.0638961038961)}] plot[domain=4.260865552203142:5.11890326389395,variable=\t]({1.*4.003061875056446*cos(\t r)+0.*4.003061875056446*sin(\t r)},{0.*4.003061875056446*cos(\t r)+1.*4.003061875056446*sin(\t r)});
\end{tikzpicture}
\end{minipage}

The Hopf map $S^3\to S^2$ takes $(z_1,z_2)\in S^3\subseteq \BC\times \BC$ to $\frac{z_1}{z_2}\in \BC\cup\{\infty\}=S^2$.
Associated to a fibration there is the Serre spectral sequence \ref{ss:Serrehomology}. In this case, its second page consists of the homology of $S^2$ with coefficients in the homology of $S^1$, and it \emph{converges} to the homology of $S^3$:
\[
E^2_{p,q}=H_p(S^2;H_q(S^1;\BZ))\Rightarrow H_{p+q}(S^3;\BZ).
\]
So we have:
\[
E^2_{p,q}=H_p(S^2;H_q(S^1;\BZ))=\begin{cases} \BZ\text{, if $p=0,2$ and $q=0,1,$}\\0\text{, otherwise,}\end{cases}
\]
and $E^2$ looks like this:
\[
\xymatrix@=0pt{
3 & 0   &  0 & 0   & 0\\
2 & 0   &  0 & 0   & 0\\
1 & \BZ &  0 & \BZ & 0\\
0\ar@{-}[]+D;[rrrr]+DR & \BZ &  0 &  \BZ & 0\\
{\phantom{0}}\ar@{-}[]+R;[uuuu]+UR  &  0  &  1 & 2   & 3
}
\]
The only possible non-trivial differential $d_2$ in $E^2$ is $d_2\colon E^2_{2,0}=\BZ\to E^2_{0,1}=\BZ$. All higher differentials $d_3,d_4,\ldots$ must be zero because of the location of the $\BZ$'s in the diagram. This ensures that point \eqref{enumerate:convergencestabilizes} above holds. As $d_2\colon \BZ\to \BZ$ is a $\BZ$-homomorphism, it is determined by $d_2(1)$. Depending on this value, we get the following $E^\infty$-page:

\begin{minipage}{0.3\textwidth}
\begin{align*}
&\xymatrix@=0pt{
3 & 0   &  0 & 0   & 0\\
2 & 0   &  0 & 0   & 0\\
1 & \BZ &  0 & \BZ & 0\\
0\ar@{-}[]+D;[rrrr]+DR & \BZ &  0 & \BZ & 0\\
{\phantom{0}}\ar@{-}[]+R;[uuuu]+UR  &  0  &  1 & 2   & 3
}\\
&\hspace{10pt}d_2(1)=0
\end{align*}
\end{minipage}
\begin{minipage}{0.3\textwidth}
\begin{align*}
&\xymatrix@=0pt{
3 & 0   &  0 & 0   & 0\\
2 & 0   &  0 & 0   & 0\\
1 & 0 &  0 & \BZ & 0\\
0\ar@{-}[]+D;[rrrr]+DR & \BZ &  0 & 0 & 0\\
{\phantom{0}}\ar@{-}[]+R;[uuuu]+UR  &  0  &  1 & 2   & 3
}\\
&\hspace{10pt}d_2(1)=\pm 1
\end{align*}
\end{minipage}
\begin{minipage}{0.3\textwidth}
\begin{align*}
&\xymatrix@=0pt{
3 & 0   &  0 & 0   & 0\\
2 & 0   &  0 & 0   & 0\\
1 & \BZ_n &  0 & \BZ & 0\\
0\ar@{-}[]+D;[rrrr]+DR & \BZ &  0 & 0 & 0\\
{\phantom{0}}\ar@{-}[]+R;[uuuu]+UR  &  0  &  1 & 2   & 3
}\\
&\hspace{10pt}d_2(1)=n,n\neq 0,\pm 1
\end{align*}
\end{minipage}

Now, point \eqref{enumerate:convergencefiniteextensions} above means that each $\BZ$-module $H_n(S^3;\BZ)$ can be recovered from the ``diagonal'' $\BZ$-modules $\{E^\infty_{p,q}\}_{p+q=n}$:
\[
\xymatrix@=0pt{
3 & 0   &  0 & 0   & 0\\
2 & 0   &  0 & 0   & 0\\
1 & ? &  0 & \BZ & 0\\
0\ar@{-}[]+D;[rrrr]+DR & \ar@{~>}[]+L;[rdd]+U \BZ &  \ar@{~>}[lu]+L;[rdd]+U 0 & 
\ar@{~>}[uull]+L;[rdd]+U {?} & \ar@{~>}[uuulll]+L;[rdd]+U 0\\
{\phantom{0}}\ar@{-}[]+R;[uuuu]+UR  &  0  &  1 & 2   & 3\\
&&H_0(S^3;\BZ)&H_1(S^3;\BZ)&H_2(S^3;\BZ)&H_3(S^3;\BZ)
}
\]
%More precisely, there exist abelian groups $A,B,C$ and short exact sequences:
%\begin{align*}
%0\to &H_0(S^3;\BZ)\to E^\infty_{0,0},\\
%E^\infty_{0,1}\to &H_1(S^3;\BZ)\to E^\infty_{1,0},\\
%A\to &H_2(S^3;\BZ)\to E^\infty_{2,0}\text{, }E^\infty_{0,2}\to A\to E^\infty_{1,1},\\
%B\to &H_3(S^3;\BZ)\to E^\infty_{3,0}\text{, }C\to B\to E^\infty_{2,1}\text{, }E^\infty_{0,3}\to C\to E^\infty_{1,2}.
%\end{align*}
In the three different cases discussed above, $d_2(1)=0$, $d_2(1)=\pm 1$, $d_2(1)=n,n\neq 0,\pm 1$, we would get the following values:

\begin{minipage}{0.3\textwidth}
\begin{align*}
&H_0(S^3;\BZ)=\BZ\\
&H_1(S^3;\BZ)=\BZ\\
&H_2(S^3;\BZ)=\BZ\\
&H_3(S^3;\BZ)=\BZ\\
&(d_2(1)=0)
\end{align*}
\end{minipage}
\begin{minipage}{0.3\textwidth}
\begin{align*}
&H_0(S^3;\BZ)=\BZ\\
&H_1(S^3;\BZ)=0\\
&H_2(S^3;\BZ)=0\\
&H_3(S^3;\BZ)=\BZ\\
&(d_2(1)=\pm 1)
\end{align*}
\end{minipage}
\begin{minipage}{0.3\textwidth}
\begin{align*}
&H_0(S^3;\BZ)=\BZ\\
&H_1(S^3;\BZ)=\BZ_n\\
&H_2(S^3;\BZ)=0\\
&H_3(S^3;\BZ)=\BZ\\
&(d_2(1)=n,n\neq 0,\pm 1).
\end{align*}
\end{minipage}

As we know the homology of $S^3$, we deduce that we must have $d_2(1)=\pm 1$ and that $E^\infty$ contains exactly two entries different from zero.
\end{exm}

\begin{rmk}
Note that the double indexing in examples \ref{example:semisimplicialmodelofS^2} and \ref{example:HopffibrationS1S3S2} are different. The one used in the latter case is customary. 
\end{rmk}

\begin{exe}
Deduce using the same arguments that if $S^l\to S^m\to S^n$ is a fibration then $n=l+1$ and $m=2n-1$.
\end{exe}

Next we unravel what point \eqref{enumerate:convergencefiniteextensions} above means with more detail: You can recover the $R$-module $H_n$ from the ``diagonal'' $R$-modules $\{E^\infty_{p,q}\}_{p+q=n}$ via a finite number of extensions.
\[
\xymatrix@=0pt{
 \ar@{~>}[]+LU;[ddddrrrr]+RD\bullet & \bullet & \bullet & \bullet & \bullet \\
 \bullet & \bullet & \bullet & \bullet & \bullet \\
 \bullet & \bullet &\bullet & \bullet & \bullet \\
 \bullet & \bullet & \bullet & \bullet & \bullet \\
 \bullet & \bullet & \bullet & \bullet & \bullet \\
         &         &         &         &         & H_n
}
\]
More precisely, there exist numbers $s\leq r$ and a finite \emph{increasing} filtration of $H_n$ by $R$-modules,
\[
0=A_{s-1}\subseteq A_s\subseteq A_{s+1}\subseteq\ldots\subseteq A_{r-1}\subseteq A_r=H_n,
\]
together with short exacts sequences of $R$-modules:
\begin{align}
&A_{r-1}\to H_n \to E^\infty_{r,n-r}\label{equ:filtrationhoonetotaldimension}\\
&A_{r-2}\to A_{r-1} \to E^\infty_{r-1,n-r+1}\nonumber\\
&\ldots\nonumber\\
&A_s\to A_{s+1}\to E^\infty_{s+1,n-s-1}\nonumber\\
&0\to A_s\to E^\infty_{s,n-s}.\nonumber
\end{align}
So, starting from the bottom, $A_s=E^\infty_{s,n-s}$, one expects to find out from these extensions what the $R$-modules  $A_{s+1},\ldots,A_{r-1},A_r=H_n$ are. This is in general not possible without further information. For instance, in Example \ref{example:HopffibrationS1S3S2} we could deduce the homology groups $H_*$ because there was only a non-zero entry on each diagonal of $E^\infty$. In the particular case of $R$ being a field, one can also deduce all terms because all what is needed is their dimension. The next example exhibits the kind of extension problems one can find when $R$ is not a field.

\begin{rmk}
If $H$ is graded and $a\in H_p$ is homogeneous then $|a|=p$ is called its degree or its total degree. If $E$ is bigraded and $a\in E_{p,q}$ then $|a|=(p,q)$ is called its bidegree and $|a|=p+q$ is called its total degree. So Equation \eqref{equ:filtrationhoonetotaldimension} says that elements from $E^{\infty}_{*,*}$ of a given total degree contribute to $H_*$ on the same total degree.
\end{rmk}

\begin{rmk} Already in the previous examples \ref{example:semisimplicialmodelofS^2} and \ref{example:HopffibrationS1S3S2}, plenty of the slang used by ``spectral sequencers'' becomes useful. For instance, regarding Example \ref{example:semisimplicialmodelofS^2}, one says that ``$B$ dies killing $c$'' and that ``$c$ dies killed by $B$'', being the reason that some differential ($d_0$) applied to $B$ is exactly $c$. Also, as $d_n=0$ for all $n\geq 1$, one says that the spectral sequence ``collapses'' at $E^1$. The terminology refers to the fact that $E^n=E^1$ for all $n\geq 2$. In Example \ref{example:HopffibrationS1S3S2}, the $\BZ$ in $E^2_{2,0}$ dies killing the $\BZ$ in $E^2_{0,1}$ and the spectral sequence collapses at $E^3$.
\end{rmk}

\begin{exm}\label{example:C2C4C2F2homol}
The homology of the cyclic group of $2^i$ elements, $C_{2^i}$, with trivial coefficients $\BZ$ or $\BZ_2$ is given by
\[
H_n(C_{2^i};\BZ)=\begin{cases}\BZ\text{, $n=0$,}\\ \BZ_{2^i}\text{, $n$ odd,}\\0\text{, otherwise}
\end{cases}\text{ and }
H_n(C_{2^i};\BZ_2)=\BZ_2\text{ for all $n\geq 0$}.
\]
The groups $C_2$ and $C_4$ fit into the short exact sequence 
\[
C_2\to C_4\to C_2.
\]
Associated to a short exact sequence of groups there is the Lyndon-Hochschild-Serre spectral sequence \ref{ss:LHShomology}. In this case, its second page consists of the homology of $C_2$ with coefficients in the homology of $C_2$, and it \emph{converges} to the homology of $C_4$:
\[
E^2_{p,q}=H_p(C_2;H_q(C_2;\BZ))\Rightarrow H_{p+q}(C_4;\BZ).
\]
By the first paragraph, $E^2$ looks like this:
\[
\xymatrix@=0pt{
{\phantom{3}} & \BZ_2  & \BZ_2& \BZ_2   & \BZ_2& \BZ_2   & \BZ_2\\
{\phantom{2}} & 0   &  0 & 0   & 0 & 0 &0\\
{\phantom{1}} & \BZ_2  & \BZ_2& \BZ_2   & \BZ_2& \BZ_2   & \BZ_2\\
{\phantom{0}} & \BZ &  \BZ_2  &  0 & \BZ_2 &  0 & \BZ_2\\
{\phantom{0}}\ar@{-}[]+U;[rrrrrr]+UR\ar@{-}[]+R;[uuuu]+UR  &  {\phantom{0}}  & {\phantom{1}} & {\phantom{2}}   & {\phantom{3}}&{\phantom{4}} &{\phantom{5}}
}
\]
Now, we know that $H_n(C_4;\BZ)=0$ for $n>0$ even. Then the $\BZ_2$'s in $E^2_{1,1}$, $E^2_{3,1}$, $E^2_{5,1}$, etc, must die and their only chance is being killed by $d_2$ from  the $\BZ_2$'s in $E^2_{3,0}$, $E^2_{5,0}$, $E^2_{7,0}$, etc: 
\[
\xymatrix@=0pt{
{\phantom{3}} & \BZ_2  & \BZ_2& \BZ_2   & \BZ_2& \BZ_2   & \BZ_2\\
{\phantom{2}} & 0   &  0 & 0   & 0 & 0 &0 \\
{\phantom{1}} & \BZ_2  & \BZ_2& \BZ_2   & \BZ_2& \BZ_2   & \BZ_2\\
{\phantom{0}} & \BZ &  \BZ_2  &  0 & \ar[]+LC;[llu]+DR \BZ_2 &  0 & \ar[]+LC;[llu]+DR\BZ_2\\
{\phantom{0}}\ar@{-}[]+U;[rrrrrr]+UR\ar@{-}[]+R;[uuuu]+UR  &  {\phantom{0}}  & {\phantom{1}} & {\phantom{2}}   & {\phantom{3}}&{\phantom{4}} &{\phantom{5}} 
}
\]
There are no other possible non-trivial differentials $d_2$ and hence $E^3$ is as follows:
\[
\xymatrix@=0pt{
{\phantom{5}} & \BZ_2  & \BZ_2& \BZ_2   & \BZ_2& \BZ_2   & \BZ_2& \BZ_2\\
{\phantom{4}} & 0   &  0 & 0   & 0 & 0 &0&0\\
{\phantom{3}} & \BZ_2  & \BZ_2& \BZ_2   & \BZ_2& \BZ_2   & \BZ_2& \BZ_2\\
{\phantom{2}} & 0   &  0 & 0   & 0 & 0 &0&0\\
{\phantom{1}} & \BZ_2  & 0& \BZ_2   & 0& \BZ_2   & 0&\BZ_2\\
{\phantom{0}} & \BZ &  \BZ_2  &  0 & 0 &  0 & 0&0\\
{\phantom{0}}\ar@{-}[]+U;[rrrrrrr]+UR\ar@{-}[]+R;[uuuuuu]+UR  &  {\phantom{0}}  & {\phantom{1}} & {\phantom{2}}   & {\phantom{3}}&{\phantom{4}} &{\phantom{5}}&{\phantom{6}}
}
\]
Again there are $\BZ_2$'s in diagonals contributing to $H_n(C_4;\BZ)$ with $n>0$ even. A careful analysis shows that all must die killed by $d_3$:
\[
\xymatrix@=0pt{
{\phantom{5}} & \BZ_2  & \BZ_2& \BZ_2   & \BZ_2& \BZ_2   & \BZ_2&\BZ_2\\
{\phantom{4}} & 0   &  0 & 0   & 0 & 0 &0&0\\
{\phantom{3}} & \BZ_2  & \BZ_2& \BZ_2   & \BZ_2&\ar[]+LC;[llluu]+DR \BZ_2   & \BZ_2&\ar[]+LC;[llluu]+DR\BZ_2\\
{\phantom{2}} & 0   &  0 & 0   & 0 & 0 &0&0\\
{\phantom{1}} & \BZ_2  & 0& \BZ_2   & 0& \ar[]+LC;[llluu]+DR\BZ_2   &0&\ar[]+LC;[llluu]+DR\BZ_2\\
{\phantom{0}} & \BZ &  \BZ_2  &  0 & 0 &  0 & 0&0\\
{\phantom{0}}\ar@{-}[]+U;[rrrrrrr]+UR\ar@{-}[]+R;[uuuuuu]+UR  &  {\phantom{0}}  & {\phantom{1}} & {\phantom{2}}   & {\phantom{3}}&{\phantom{4}} &{\phantom{5}} &{\phantom{6}}
}
\]
Hence, $E^4$ looks as follows. 
\[
\xymatrix@=0pt{
{\phantom{3}} & \BZ_2  &0& \BZ_2   & 0& 0   & 0&0\\
{\phantom{2}} & 0   &  0 & 0   & 0 & 0 &0&0\\
{\phantom{3}} & \BZ_2  &0& \BZ_2   & 0& 0   & 0&0\\
{\phantom{2}} & 0   &  0 & 0   & 0 & 0 &0&0\\
{\phantom{1}} & \BZ_2  & 0& \BZ_2   & 0&0   & 0&0\\
{\phantom{0}} & \BZ &  \BZ_2  &  0 &0 &  0 &0&0\\
{\phantom{0}}\ar@{-}[]+U;[rrrrrrr]+UR\ar@{-}[]+R;[uuuuuu]+UR  &  {\phantom{0}}  & {\phantom{1}} & {\phantom{2}}   & {\phantom{3}}&{\phantom{4}} &{\phantom{5}}&{\phantom{6}}
}
\]
It is clear from the positions of the non-zero entries in $E^4$ that the rest of the differentials $d_4,d_5,d_6,\ldots$ must be zero. So the spectral sequence collapses at $E^4=E^\infty$. A posteriori, from $E^\infty$, one can deduce that $H_0(C_4;\BZ)=\BZ$ and that $H_n(C_4;\BZ)=0$ for $n$ even, $n>0$. For $n$ odd, we have an extension of $\BZ$-modules:
\[
\BZ_2\to H_n(C_4;\BZ)\to \BZ_2.
\]
There are two solutions, either $H_n(C_4;\BZ)=\BZ_4$ or $H_n(C_4;\BZ)=\BZ_2\times \BZ_2$, and one would need extra information to decide which one is the right one.
\end{exm}

\begin{exe}\label{exe:loopsonS^n}
From the fibration $\Omega S^n\to *\to S^n$ and the Serre spectral sequence \ref{ss:Serrehomology}, deduce the integral homology of the loops on the sphere $\Omega S^n$. Solution in \ref{exe:loopsonS^nsolution}.
\end{exe}

\section{Cohomological type spectral sequences}
\label{ss4:cohotype}

A cohomological type spectral sequences consists of a sequence of bigraded differential modules  $\{E_r,d_r\}_{r\geq 2}$ such that $E_{r+1}$ is the cohomology of $(E_r,d_r)$ and such that $d_r$ has bidegree $(r,1-r)$, $d_r\colon E^r_{p,q}\to E^r_{p+r,q+1-r}$, as in these pictures:

\begin{minipage}{0.3\textwidth}
\begin{align*}
&\xymatrix@=0pt{
 \bullet & \bullet & \bullet & \bullet & \bullet \\
 \bullet & \bullet & \bullet & \bullet & \bullet \\
 \bullet & \bullet & \ar@{<-}[]+C;[llu]+C\bullet & \bullet & \bullet \\
 \bullet & \bullet & \bullet & \bullet & \bullet \\
 \bullet & \bullet & \bullet & \bullet & \bullet \\
}\\
&\hspace{10pt}(E_2,d_2)
\end{align*}
\end{minipage}
\begin{minipage}{0.3\textwidth}
\begin{align*}
&\xymatrix@=0pt{
 \bullet & \bullet & \bullet & \bullet & \bullet \\
 \bullet & \bullet & \bullet & \bullet & \bullet \\
 \bullet & \bullet & \bullet & \bullet & \bullet \\
 \bullet & \bullet & \bullet & \ar@{<-}[]+C;[llluu]+C\bullet & \bullet \\
 \bullet & \bullet & \bullet & \bullet & \bullet \\
}\\
&\hspace{10pt}(E_3,d_3)
\end{align*}
\end{minipage}
\begin{minipage}{0.3\textwidth}
\begin{align*}
&\xymatrix@=0pt{
 \bullet & \bullet & \bullet & \bullet & \bullet \\
 \bullet & \bullet & \bullet & \bullet & \bullet \\
 \bullet & \bullet & \bullet & \bullet & \bullet \\
 \bullet & \bullet & \bullet & \bullet &  \ar@{<-}[]+C;[lllluuu]+C\bullet \\
 \bullet & \bullet & \bullet & \bullet & \bullet \\
}\\
&\hspace{10pt}(E_4,d_4)
\end{align*}
\end{minipage}

In this situation, $H=H^*$ is the cohomology of some  \emph{cochain} complex, $C^*$, and convergence is written $E_2^{p,q}\Rightarrow H^{p+q}$. As for homological type, convergence means that the pages eventually stabilize at each particular position ($p,q$) and that  you can recover the $R$-module $H^n$ from the ``diagonal'' $R$-modules $\{E_\infty^{p,q}\}_{p+q=n}$ via a finite number of extensions. More precisely, there exist numbers $s\leq r$ and a finite \emph{decreasing} filtration of $H^n$ by $R$-modules,
\[
0=A^{r+1}\subseteq A^r\subseteq A^{r-1}\subseteq\ldots\subseteq A^{s+1}\subseteq A^s=H^n,
\]
together with short exact sequences of $R$-modules:
\begin{align}
&A^{s+1}\to H^n \to E_\infty^{s,n-s}\label{equ:filtrationcohoonetotaldimension}\\
&A^{s+2}\to A^{s+1} \to E_\infty^{s+1,n-s-1}\nonumber\\
&\ldots\nonumber \\
&A^r\to A^{r-1}\to E_\infty^{r-1,n-r+1}\nonumber \\
&0\to A^r\to E_\infty^{r,n-r}.\nonumber
\end{align}
\begin{rmk}\label{rmk:upsidedownindexing}
Pay attention to the difference between \eqref{equ:filtrationhoonetotaldimension}  and
 \eqref{equ:filtrationcohoonetotaldimension}: the indexing is upside down.
\end{rmk}

\begin{exm}\label{example:KtheoryCP^n}
The functor complex $K$-theory satisfies that $K^n(X)=K^{n+2}(X)$ for any compact topological space $X$. For a point, we have $K^0(*)=\BZ$ and $K^1(*)=0$. Let $\BC P^n$ be the complex projective space of dimension $n$ and consider the fibration
\[
*\to \BC P^n \stackrel{id}\longrightarrow \BC P^n.
\]
Associated to a fibration there is the Atiyah-Hirzebruch spectral sequence \ref{ss:AtiyaHirzebruchcohomology}. In this case, the second page is given by the cohomology of $\BC P^n$ with coefficients in the $K$-theory of a point, and it converges to the $K$-theory of $\BC P^n$:
\[
E_2^{p,q}=H^p(\BC P^n; K^q(*))\Rightarrow K^{p+q}(\BC P^n).
\] 
The integral cohomology groups of $\BC P^n$ are given by $H^p(\BC P^n;\BZ)=\BZ$ for $0\leq p\leq n$, $p$ even, and $0$ otherwise. So $E_2^{p,q}=\BZ$ if $p$ and $q$ are even and $0\leq p\leq 2n$, and it is $0$ otherwise. This is the picture for $\BC P^3$:
\[
\xymatrix@=0pt{
{\phantom{4}} & \BZ &  0  & \BZ & 0 & \BZ & 0 &\BZ \\
{\phantom{3}} & 0  & 0& 0   & 0& 0   & 0&0\\
{\phantom{2}}& \BZ &  0  & \BZ & 0 & \BZ & 0&\BZ \\
{\phantom{1}}\ar@{-}[]+U;[rrrrrrr]+UR & 0  & 0& 0   & 0& 0   & 0&0\\
{\phantom{0}} & \BZ &  0  & \BZ & 0 & \BZ & 0 &\BZ\\
{\phantom{0}}\ar@{-}[]+R;[uuuuu]+UR  &  0  & 0 & 0   & 0&0 &0 &0
 }
\]
Because for every $r$, either $r$ or $r-1$ is odd, it follows that $d_r=0$ for all $r\geq 2$. Hence the spectral sequence collapses at $E_2=E_\infty$. Because $\BZ$ is a free abelian group, all extension problems have a unique solution and we deduce that:
\[
K^k(\BC P^n)=\begin{cases}\BZ^{n+1} \text{, if $k$ is even,}\\0\text{, otherwise.}\end{cases}
\]
\end{exm}

\section{Additional structure: algebra}
\label{ss5:algebra}

In some cases, the cohomology groups we want to compute, $H$, have some additional structure. For instance, there may be a product
\[
H\otimes H\to H
\]
that makes $H$ into an algebra. A \emph{spectral sequences of algebras} is a spectral sequence $\{E_r,d_r\}_{r\geq 2}$ such that $E_r$ is a differential bigraded algebra, i.e., there is a product,
\[
E_r\otimes E_r\to E_r,
\]
and the product in $E_{r+1}$ is that induced by the product in $E_r$ after taking cohomology with respect to $d_r$. A spectral sequence of algebras \emph{converge as an algebra} to the algebra $H$ if it converges as a spectral sequence and there is a decreasing filtration of $H^*$,
\[
\ldots \subseteq F^{n+1}\subseteq F^n \subseteq F^{n-1}\subseteq \ldots \subseteq H^*,
\]
which is compatible with the product in $H$,
\[
F^n\cdot F^m\subseteq F^{n+m},
\]
and which satisfies the conditions below. When one restricts this filtration to each dimension, one gets the filtration explained in \eqref{equ:filtrationcohoonetotaldimension}. More precisely, if we define $F^pH^{p+q}=F^p\cap H^{p+q}$, we have short exact sequences 
\[
F^{p+1} H^{p+q}\to F^p H^{p+q}\to E_\infty^{p,q}.
\]
Then there are two products in the page $E_\infty$:
\begin{enumerate}[(a)]
\setcounter{enumi}{2}
\item \label{enumerate:Einftyproductislimit}
The one induced by $E_\infty$ being the limit of the algebras $E_2,E_3,\ldots$,
\item \label{enumerate:Einftyproductisbyfiltration}
The one induced by the filtration as follows: For elements $\overline a$ and $\overline b$ belonging to $E_\infty^{p,q}=F^p H^{p+q}/F^{p+1} H^{p+q}$ and $E_\infty^{p',q'}=F^{p'} H^{p'+q'}/F^{p'+1} H^{p'+q'}$ respectively, we set:
\[
\overline a\cdot \overline b=\overline{a\cdot b}\in E_\infty^{p+p',q+q'}=F^{p+p'} H^{p+p'+q+q'}/F^{p+p'+1} H^{p+p'+q+q'}.
\]
\end{enumerate}
We say that the spectral sequence \emph{converge as an algebra} if these two products are equal.
\begin{rmk} 
Note that, by definition, the graded and bigraded products satisfy $H^n\cdot H^m\subseteq H^{n+m}$ and $E^{p,q}_r\cdot E^{p',q'}_r\subseteq E^{p+p',q+q'}_r$ for $r\in \BN \cup \{\infty\}$. All algebras we consider are graded (or bigraded) commutative, i.e.,
\begin{equation}\label{equ:gradedalgebrarelations}
ab=(-1)^{|a||b|}ba,
\end{equation}
where $|a|$ and $|b|$ are the degrees (total degrees) of the homogeneous elements $a$ and $b$ respectively. Moreover, every differential $d$ is a derivation, i.e.,
\[
d(ab)=d(a)b+(-1)^{|a|}ad(b).
\]
Finally, if $C^{*,*}$ is a bigraded module or a double cochain complex, by $\total(C)$ we denote the graded algebra or the cochain complex such that 
\[
\total(C)^n=\bigoplus_{p+q=n} C^{p,q}.
\]
\end{rmk}
\begin{exm}\label{example:C2C4C2F2}
Consider the short exact sequence of cyclic groups:
\[
C_2\to C_4\to C_2.
\]
The cohomology ring of $C_2$ with coefficients in the field of two elements is given by $H^*(C_2;\BF_2)=\BF_2[x]$ with $|x|=1$. For the extension above there is the Lyndon-Hochschild-Serre spectral sequence of algebras converging as an algebra \ref{ss:LHScohomology},
\[
E_2^{p,q}=H^p(C_2,H^q(C_2;\BF_2))\Rightarrow H^{p+q}(C_4;\BF_2).
\]
As the extension is central, the corner of the page $E_2=\BF_2[x,y]$ has the following generators:
\[
\xymatrix@=0pt{
{\phantom{3}} & y^3  & y^3x & y^3x^2 & y^3x^3\\
{\phantom{2}} & y^2  & y^2x & y^2x^2 & y^2x^3\\
{\phantom{1}} & y  & yx & yx^2 & yx^3\\
{\phantom{0}} & 1 &  x  &  x^2 & x^3 \\
{\phantom{0}}\ar@{-}[]+U;[rrrr]+UR\ar@{-}[]+R;[uuuu]+UR  &  {\phantom{0}}  & {\phantom{1}} & {\phantom{2}}   & {\phantom{3}}&{\phantom{4}} &{\phantom{5}}
}
\]
Because the extension is non-split, we know that $d_2(y)=x^2$. Another way of deducing this is to use the fact $H^1(C_4;\BF_2)=\BF_2$. Then, as the terms $E_\infty^{1,0}=\langle x\rangle$ and $E_\infty^{0,1}$ contribute to $H^1(C_4;\BF_2)$, the $y$ must die, and $d_2(y)=x^2$ is its only chance. From here, we can deduce the rest of the differentials as $d_2$ is a derivation, for instance:
\begin{align*}
&d_2(yx)=d_2(y)x+yd_2(x)=x^3\text{, as $d_2(x)=0$,}\\
&d_2(y^2)=d_2(y)y+yd_2(y)=2d_2(y)y=0\text{, as we work over $\BF_2$.}
\end{align*}
Analogously, one deduces that $d_2(y^{2i+1}x^j)=y^{2i}x^{j+2}$ and that $d_2(y^{2i}x^j)=0$. Hence $d_2$ kills elements as follows,
\[
\xymatrix@=0pt{
{\phantom{3}} & y^3\ar[]+C;[rrd]+C  & y^3x\ar[]+C;[rrd]+C & y^3x^2 & y^3x^3\\
{\phantom{2}} & y^2  & y^2x & y^2x^2 & y^2x^3\\
{\phantom{1}} & y\ar[]+C;[rrd]+C  & yx \ar[]+C;[rrd]+C& yx^2 & yx^3\\
{\phantom{0}} & 1 &  x  &  x^2 & x^3 \\
{\phantom{0}}\ar@{-}[]+U;[rrrr]+UR\ar@{-}[]+R;[uuuu]+UR  &  {\phantom{0}}  & {\phantom{1}} & {\phantom{2}}   & {\phantom{3}}&{\phantom{4}} &{\phantom{5}}
}
\]
and $E_3$ is
\[
\xymatrix@=0pt{
{\phantom{3}} & 0& 0 & 0 & 0\\
{\phantom{2}} & y^2  & y^2x & 0 & 0\\
{\phantom{1}} & 0  & 0& 0 & 0\\
{\phantom{0}} & 1 &  x  &  0& 0 \\
{\phantom{0}}\ar@{-}[]+U;[rrrr]+UR\ar@{-}[]+R;[uuuu]+UR  &  {\phantom{0}}  & {\phantom{1}} & {\phantom{2}}   & {\phantom{3}}&{\phantom{4}} &{\phantom{5}}
}
\]
By degree reasons, there cannot be any other non-trivial differentials and hence the spectral sequence collapses at $E_3=E_\infty$. Moreover, the bigraded algebra $E_\infty$ is given by $E_{\infty}=\BF_2[x,x']/(x^2)=\Lambda(x)\otimes \BF_2[x']$, where $x'=y^2$, $|x|=(1,0)$ and $|x'|=(0,2)$. In fact, $H^*(C_4;\BF_2)=\Lambda(z)\otimes \BF_2[z']$ with $|z|=1$, $|z'|=2$.
\end{exm}

\begin{exe}
Compare last example and Example \ref{example:C2C4C2F2homol}.
\end{exe}

\begin{exe}\label{exe:C3C9C3coho}
Do the analogous computation for the extension $C_3\to C_9\to C_3$ given that $H^*(C_3;\BF_3)=\Lambda(y)\otimes \BF_3[x]$ with $|y|=1$, $|x|=2$. Solution in \ref{exe:C3C9C3cohosolution}.
\end{exe}

In general we cannot recover the graded algebra structure in $H^*$ from the bigraded algebra structure of $E_\infty$.  This problem is known as the \emph{lifting problem}. The next theorem gives some special conditions under we can reconstruct $H^*$ from $E_\infty^{*,*}$. They apply to Exercise \ref{exe:C3C9C3coho}.

\begin{thm}[{\cite[Example 1.K, p. 25]{McCleary}}]\label{thm:freebigradedfreegraded}
If $E_\infty$ is a free, graded-commutative, bigraded algebra, then $H^*$ is
a free, graded commutative algebra isomorphic to $\total(E_\infty)$.
\end{thm}

A free graded (or bigraded) commutative algebra is the quotient of the free algebra on some graded (bigraded) symbols modulo the relations \eqref{equ:gradedalgebrarelations}. Theorem \ref{thm:freebigradedfreegraded} translates as that, if $x_1,\ldots,x_r$ are the free generators of $E_\infty$ with bidegrees $(p_1,q_1),\ldots,(p_r,q_r)$, then $H^*$ is  a free graded commutative algebra on generators $z_1,\ldots,z_r$ of degrees $p_1+q_1,\ldots,p_r+q_r$. The next example is the integral version of the earlier Example \ref{example:C2C4C2F2} and here the lifting problem becomes apparent.

\begin{exm}\label{example:C2C4C2Z}.
Consider the short exact sequence of cyclic groups:
\[
C_2\to C_4\to C_2.
\]
The integral cohomology ring of the cyclic group $C_2$ is given by $H^*(C_2;\BZ)=\BZ[x]/(2x)$, with $|x|=2$. The Lyndon-Hochschild-Serre spectral \ref{ss:LHScohomology} is
\[
E_2^{p,q}=H^p(C_2,H^q(C_2;\BZ))\Rightarrow H^{p+q}(C_4;\BZ).
\]
As the extension is central, the corner of the page $E_2=\BF_2[x,y]$ has the following generators:
\[
\xymatrix@=0pt{
{\phantom{2}} & y^2 &0& y^2x &0& y^2x^2 &0& y^2x^3\\
{\phantom{1}} & 0   &0& 0    &0&0       &0& 0\\
{\phantom{1}} & y   &0& yx   &0& yx^2   &0& yx^3\\
{\phantom{1}} & 0   &0& 0    &0&0       &0& 0\\
{\phantom{0}} & 1   &0&  x   &0&  x^2   &0& x^3 \\
{\phantom{0}}\ar@{-}[]+U;[rrrrrrr]+UR\ar@{-}[]+R;[uuuuu]+UR  &  {\phantom{0}}  & {\phantom{1}} & {\phantom{2}}   & {\phantom{3}}&{\phantom{4}} &{\phantom{5}}&{\phantom{6}} &{\phantom{7}}
}
\]
Recall that $|x|=(2,0)$ and $|y|=(0,2)$. As in Example \ref{example:KtheoryCP^n}, parity implies that the spectral sequence collapses at $E_2=E_\infty$ and hence  $E_\infty=\BZ[x,y]/(2x,2y)$. Nevertheless, $H^*(C_4;\BZ)=\BZ[z]/(4z)$ with $|z|=2$. 
%A posteriori, we see that the solution to the extension problem \ref{equ:filtrationcohoonetotaldimension} for total degree $2$,
%\[
%\BZ_2=\langle x\rangle \to H^2(C_4;\BZ)\cong \BZ_4=\langle z\rangle\to \BZ_2=\langle y\rangle,
%\]
%is to take $x=2z$ and $y=\bar z$.
\end{exm}

In the next section, we deepen into the lifting problem to pin down the relation between the bigraded algebra $E_\infty$ and the graded algebra $H^*$. We finish this section with a topological example.

\begin{exm}\label{example:cohoK(Z,n)}
We are about to determine the ring $H^*(K(\BZ,n);\BQ)$, where $K(\BZ,n)$ is an Eilenberg-MacLane space, i.e., their homotopy groups satisfy
\[
\pi_i(K(\BZ,n))=\begin{cases}\BZ\text{, if $i=n$,}\\0\text{, otherwise.}\end{cases}
\]
So, we know that $K(\BZ,1)=S^1$ and hence $H^*(K(\BZ,1);\BQ)=\Lambda(z)$ with $|z|=1$. We prove by induction that 
\[
H^*(K(\BZ,n);\BQ)=\begin{cases} \BQ[z]\text{, if $n$ is even,}\\\Lambda(z)\text{, if $n$ is odd,}\end{cases}
\]
where $|z|=n$. We consider the fibration $K(\BZ,n-1)\to * \to K(\BZ,n)$ and its Serre spectral sequence \ref{ss:Serrecohomology}. Notice that $K(\BZ,n)$ is simply connected as $n\geq 2$. Then $E_2^{p,q}=H^p(K(\BZ,n);\BQ)\otimes H^q(K(\BZ,n-1);\BQ)$. As $E_2^{*,*}\Rightarrow H^*(*;\BQ)=\BQ$, all terms but the $\BZ$ in $(0,0)$ must disappear. If $n$ is even then $H^*(K(\BZ,n-1);\BQ)=\Lambda(z)$ with $|z|=n-1$. Moreover, the only chance of dying for $z$ at position $(0,n-1)$ is by killing $d_n(z)=x$ at $(n,0)$. Then $zx$ at $(n,n-1)$ must die killing $d_n(zx)=d_n(z)x+(-1)^{n-1}zd_n(x)=x^2$. 
\[
\xymatrix@=0pt{
{\phantom{1}} & 0   &0& 0    &0&0       &0& 0&0&0&0\\
{\phantom{1}} & z\ar[]+DR;[rrrrdddd]+LU^{d_n}   &0& \ldots & 0 & zx \ar[]+DR;[rrrrdddd]+LU^{d_n}  &0& \ldots  &0& zx^2&0\\
{\phantom{1}} & 0   &0& \ldots & 0 &0     &0& \ldots  &0&0&0\\
{\phantom{1}} &     && \ldots    &&       & & \\
{\phantom{1}} & 0   &0& \ldots & 0 &0     &0& \ldots  &0&0&0\\
{\phantom{1}} & 1   &0& \ldots & 0 & x    &0& \ldots  &0&x^2&0\\
{\phantom{0}}\ar@{-}[]+U;[rrrrrrrrrr]+UR\ar@{-}[]+R;[uuuuuu]+UR  &  {\phantom{0}}  & {\phantom{1}} & {\phantom{2}}   & {\phantom{3}}&{\phantom{4}} &{\phantom{5}}&{\phantom{6}} &{\phantom{7}}&{\phantom{8}}&{\phantom{9}} &{\phantom{0}}
}
\]
Pushing this argument further, one may conclude that $H^*(K(\BZ,n);\BQ)=\BQ[x]$ with $|x|=n$. Now assume that $n$ is odd. Then $H^*(K(\BZ,n-1);\BQ)=\BQ[z]$ with $|z|=n-1$ and again $z$ at  $(0,n-1)$ must die killing $d_n(z)=x$ at $(n,0)$. Next, we deduce that $d_n(z^2)=d_n(z)z+(-1)^{n-1}zd_n(z)=2zx$ and hence $zx$ is killed by $\frac{1}{2}{z^2}$:
\[
\xymatrix@=0pt{
{\phantom{1}} & 0   &0& \ldots & 0 &0     &0\\
{\phantom{1}} & \frac{1}{2}z^2\ar[]+DR;[rrrrdddd]+LU^{d_n}   &0& \ldots & 0 & z^2x  &0\\
{\phantom{1}} & 0   &0& \ldots & 0 &0     &0\\
{\phantom{1}} &     && \ldots    & \\
{\phantom{1}} & 0   &0& 0    &0&0       &0&\\
{\phantom{1}} & z\ar[]+DR;[rrrrdddd]+LU^{d_n}   &0& \ldots & 0 & zx  &0\\
{\phantom{1}} & 0   &0& \ldots & 0 &0     &0\\
{\phantom{1}} &     && \ldots    & \\
{\phantom{1}} & 0   &0& \ldots & 0 &0     &0\\
{\phantom{1}} & 1   &0& \ldots & 0 & x    &0\\
{\phantom{0}}\ar@{-}[]+U;[rrrrrrrrrr]+UR\ar@{-}[]+R;[uuuuuuuuuu]+UR  &  {\phantom{0}}  & {\phantom{1}} & {\phantom{2}}   & {\phantom{3}}&{\phantom{4}} &{\phantom{5}}&{\phantom{6}} &{\phantom{7}}&{\phantom{8}}&{\phantom{9}} &{\phantom{0}}
}
\]
Following these arguments, one may conclude that $H^*(K(\BZ,n);\BQ)=\Lambda(x)$ with $|x|=n$.
\end{exm}

The next two exercise are meant to highlight the role of the base ring $R$. Computations are similar to those of Example \ref{example:cohoK(Z,n)} but with integral coefficients instead of rational coefficients.

\begin{exe}\label{exe:cohoK(Z,2)}
Determine the cohomology ring $H^*(K(\BZ,2);\BZ)$ using the fibration $K(\BZ,1)\to *\to K(\BZ,2)$. Solution in \ref{exe:cohoK(Z,2)solution}.
\end{exe}

\begin{exe}\label{exe:cohoK(Z,3)}
Determine the cohomology ring $H^*(K(\BZ,3);\BZ)$ using the fibration $K(\BZ,2)\to *\to K(\BZ,3)$ and the previous exercise. Solution in \ref{exe:cohoK(Z,3)solution}.
\end{exe}

\begin{exe}\label{exe:coholoopsS^k}
Determine the cohomology ring $H^*(\Omega S^n;\BZ)$ using the fibration $\Omega S^n\to *\to S^n$. Solution in \ref{exe:coholoopsS^ksolution}.
\end{exe}

\section{Lifting problem}
\label{ss6:lifting}

Recall the general setting for the lifting problem: we have a graded algebra $H^*$ with a filtration $F^*$ and then we consider the associated bigraded algebra $E_\infty$ with the induced product as in \eqref{enumerate:Einftyproductisbyfiltration}. So, for $F^pH^{p+q}=F^p\cap H^{p+q}$, we have $E_\infty^{p,q}\cong F^pH^{p+q}/F^{p+1}H^{p+q}$ and the product $E_\infty^{p,q}\otimes E_\infty^{p',q'}\to E_\infty^{p+p',q+q'}$ sends
\[
(a+F^{p+1} H^{p+q})\cdot (b+F^{p'+1} H^{p'+q'})=(a\cdot b)+F^{p+p'+1}H^{p+q+p'+q'}.
\]

\begin{exm}[{\cite[Example 1.J, p. 23]{McCleary}}] \label{example:twoalgebrassasmebigradedalgebra}
Consider the following two algebras:
\begin{align*}
H_1^*=\BQ[x,y,z]/(x^2,y^2,z^2,xy,xz,yz)\text{ with $|x|=7$, $|y|=8$, $|z|=15$, and}\\
H_2^*=\BQ[u,v,w]/(u^2,v^2,w^2,uv-w,uw,vw)\text{ with $|u|=7$, $|v|=8$, $|w|=15$.}
\end{align*}
The following decreasing filtrations $F_1$ and $F_2$ of $H_1$ and $H_2$ respectively are compatible with the respective products:
\begin{align*}
0=F_1^8\subset F_1^7=F_1^6=F_1^5=\langle z\rangle\subset F_1^4=F_1^3=\langle y,z\rangle\subset F_1^2=F_1^1=H_1\text{, and}\\
0=F_2^8\subset F_2^7=F_2^6=F_2^5=\langle w\rangle\subset F_2^4=F_2^3=\langle v,w\rangle\subset F_2^2\subset F_2^1=H_2.
\end{align*}
The next pictures depict the corresponding bigraded algebras $E_{\infty,1}$ and $E_{\infty,2}$:

\begin{minipage}{0.5\textwidth}
\[
\xymatrix@=0pt{
9 & 0   &  0 & 0   & 0 & 0   &  0 & 0   & 0\\
8 & 0   &  0 & 0   & 0 & 0   &  0 & \bar z   & 0\\
7 & 0   &  0 & 0   & 0 & 0   &  0 & 0   & 0\\
6 & 0   &  0 & 0   & 0 & 0   &  0 & 0   & 0\\
5 & 0   &  \bar x & 0   & 0 & 0   &  0 & 0   & 0\\
4 & 0   &  0 & 0   & \bar y & 0   &  0 & 0   & 0\\
3 & 0   &  0 & 0   & 0 & 0   &  0 & 0   & 0\\
2 & 0   &  0 & 0   & 0 & 0   &  0 & 0   & 0\\
1 & 0   &  0 & 0   & 0 & 0   &  0 & 0   & 0\\
{\phantom{0}}\ar@{-}[]+U;[rrrrrrrr]+UR\ar@{-}[]+R;[uuuuuuuuu]+UR  &  1  &  2 & 3   & 4 &  5  & 6 & 7   & 8
}
\]
\end{minipage}
\begin{minipage}{0.5\textwidth}
\[
\xymatrix@=0pt{
9 & 0   &  0 & 0   & 0 & 0   &  0 & 0   & 0\\
8 & 0   &  0 & 0   & 0 & 0   &  0 & \bar w   & 0\\
7 & 0   &  0 & 0   & 0 & 0   &  0 & 0   & 0\\
6 & 0   &  0 & 0   & 0 & 0   &  0 & 0   & 0\\
5 & 0   &  \bar u & 0   & 0 & 0   &  0 & 0   & 0\\
4 & 0   &  0 & 0   & \bar v & 0   &  0 & 0   & 0\\
3 & 0   &  0 & 0   & 0 & 0   &  0 & 0   & 0\\
2 & 0   &  0 & 0   & 0 & 0   &  0 & 0   & 0\\
1 & 0   &  0 & 0   & 0 & 0   &  0 & 0   & 0\\
{\phantom{0}}\ar@{-}[]+U;[rrrrrrrr]+UR\ar@{-}[]+R;[uuuuuuuuu]+UR  &  1  &  2 & 3   & 4 &  5  & 6 & 7   & 8
}
\]
\end{minipage}
We have  $E_{\infty,1}\cong E_{\infty,2}$ as bigraded algebras because all products of generators are zero in both cases. In the former case, this is a consequence of the products being zero in the graded algebra $H_1^*$. In the latter case, the product $\bar u\bar v$ is zero because $|\bar u|=(2,5)$, $|\bar v|=(4,4)$, $|\bar u\bar v|=(6,9)$ and $E_{\infty,2}^{6,9}=0$. What is happening is that $\bar u\bar v=0\in E_{\infty,2}^{6,9}=F_2^6H_2^{15}/F_2^7H_2^{15}$ implies that $uv\in F_2^7H_2^{15}=\langle w\rangle$. So we may deduce that $uv=\lambda w$ for some $\lambda\in\BQ$,  but we cannot know whether $\lambda =0$ or not.
\end{exm}

The following result clarifies the general situation.

\begin{thm}[{\cite[Theorem 2.1]{Carlson}}]\label{thm:Carlsonliftinggensandrels}
If $E_\infty$ is a finitely generated bigraded commutative algebra, then $H^*$ is a finitely generated graded commutative algebra. Moreover, given a presentation of $E_\infty$ by double homogeneous generators and double homogeneous relations:
\begin{enumerate}
\item Generators for $H^*$ are given by lifting the given generators of $E_\infty$, and they have the same total degree.
\item Relations  for $H^*$ are given by lifting the given relations  of $E_\infty$, and they have the same total degree.
\end{enumerate}
\end{thm}

In Example \ref{example:twoalgebrassasmebigradedalgebra}, we have $E_{\infty,1}^{*,*}=\BQ[\bar x,\bar y,\bar z]/(\bar x^2,\bar y^2,\bar z^2,\bar x\bar y,\bar x\bar z,\bar y\bar z)$ and these generators and relations lift to give $H_1^*=\BQ[x,y,z]/(x^2,y^2,z^2,xy,xz,yz)$. For $E_{\infty,2}^{*,*}=\BQ[\bar u,\bar v,\bar w]/(\bar u^2,\bar v^2,\bar w^2,\bar u\bar v,\bar u\bar w,\bar v\bar w)$,  the relation $\bar u\bar v=0$ lifts as $uv=w$. This theorem has the following immediate consequence.

\begin{cor}[{\cite[Theorem 2.1]{Carlson}}]
If $R$ is a finite field, the graded algebra structure of $H^*$ is determined by the bigraded algebra structure of $E_\infty$ within a finite number of posibilities.
\end{cor}

At this point, the original picture of the bridge and the tree should be clearer. That there is a bridge means that the spectral sequence converges, and you can cross it if you are able to figure out enough differentials. Climbing the tree implies solving extension and lifting problems. This may require extra information.

\begin{center}
\begin{picture}(100,75)(0,-20)
\qbezier(0,0)(5,15)(20,20)
%\qbezier(20,20)(27,22)(33,21)
\qbezier(33,21)(40,21)(47,19)
\qbezier(47,15)(30,10)(27,-13)
\qbezier(47,15)(53,16)(56,15)
\qbezier(47,19)(53,17)(56,15)
\qbezier(56,15)(63,10)(65,0)

\qbezier(0,0)(0,0)(27,-13)
\qbezier(38,10)(65,0)(65,0)

%RIO
\qbezier(28,-8)(28,-8)(53,-20)
%\qbezier(33,-8)(33,-8)(60,-21)
%\qbezier(40,-6)(40,-6)(67,-17)
\qbezier(38,-9)(38,-9)(65,-20)
\qbezier(40,5)(40,5)(67,-5)
\qbezier(35,6.8)(40,5)(40,5)
\qbezier(67,-5)(67,-5)(90,-13)
%\qbezier(57,-5)(57,-5)(80,-13)
\qbezier(52,-4.5)(52,-4.5)(75,-13)

\qbezier(30,-1)(30,-1)(53,-10.5)
\qbezier(35,5)(35,5)(60,-4.5)

\qbezier(-25,16)(-25,16)(1,4)
\qbezier(-20,26)(-20,26)(10,15)

\qbezier(-20,16)(-20,16)(2.7,6)
\qbezier(-10,16)(-10,16)(5,9.6)

\put(-17,-7){$E_2$}
\put(9,0){$d_2$}
\put(14,10){$d_3$}
\put(22,17){$d_4$}
\put(33,17){$\dots $}
\put(-17,35){Convergence}
\put(10,25){$\longrightarrow$}

\put(70,5){$E_\infty$}
\put(95,20){$\Big\uparrow$}
\put(105,25){extension and}
\put(105,15){lifting problems}

\qbezier(68,16)(73,25)(71,28)
\qbezier(82,16)(77,25)(79,28)

%COPA
\qbezier(60,31)(75,26)(90,31)
\qbezier(60,31)(55, 34)(60,41)
\qbezier(90,31)(95,34)(90, 41)
\qbezier(60, 41)(75,58)(90, 41)

%SOMBRA
\multiput(38.5,8)(1.5,-0.55){17}{$\cdot$}
\multiput(39.5,9)(1.5,-0.55){16}{$\cdot$}
\multiput(40.5,10)(1.5,-0.55){15}{$\cdot$}
\multiput(43,10.45)(1.5,-0.55){13}{$\cdot$}
\multiput(44,11.45)(1.5,-0.55){12}{$\cdot$}
\multiput(46.5,11.90)(1.5,-0.55){9}{$\cdot$}
\multiput(49,12.35)(1.5,-0.55){7}{$\cdot$}
\multiput(53,12.25)(1.5,-0.55){2}{$\cdot$}
\put(72,35){$H$}
\end{picture}
\end{center}

We exemplify Theorem \ref{thm:Carlsonliftinggensandrels} by means of the next example.

\begin{exm}\label{example:dihedralaswreath}
Consider the dihedral group $D_8$ described as a wreath product:
\[
D_8=C_2\wr C_2=C_2\times C_2\rtimes C_2=\langle a,b\rangle \rtimes \langle c\rangle,
\]
where $c$ interchanges $a$ and $b$: ${}^ca=b$, ${}^cb=a$. The Lyndon-Hochschild-Serre spectral sequence \ref{ss:LHScohomology} with coefficients in the field of two elements is
\[
E_2^{p,q}=H^p(C_2;H^q(C_2\times C_2;\BF_2))\Rightarrow H^{p+q}(D_8;\BF_2).
\]
\end{exm}
The extension is not central and it is not immediate to describe $E_2$. Recall that $H^*(C_2\times C_2;\BF_2) =H^*(C_2;\BF_2)\otimes H^*(C_2;\BF_2)=\BF_2[x,y]$, with $|x|=|y|=1$ . Here, $x=a^*$ and $y=b^*$ are the corresponding duals. This ring can be decomposed as a $C_2$-module into a direct sum of two type of $C_2$-modules: either the trivial $C_2$-module or the free transitive $C_2$-module:

\[
\xymatrix@=5pt{ 1\ar@{<->}@(ur,ul) & x\ar@{<->}@(ul,ur)[r]&y & x^2\ar@{<->}@(ul,ur)[r] & y^2 & xy\ar@{<->}@(ur,ul) & x^3\ar@{<->}@(ul,ur)[r] & y^3 & x^2y\ar@{<->}@(ul,ur)[r] & xy^2 &x^2y^2\ar@{<->}@(ur,ul)&\ldots}
\]
Now, the cohomology of the trivial module has been already mentioned $H^*(C_2;\BF_2)=\BF_2[z]$ with $|z|=1$. The cohomology of the free transitive module $M$ is given by the fixed points in degree $0$, $H^0(C_2;M)=M^{C_2}$, and $H^n(C_2;M)=0$ for $n>0$ (as $0\to M\stackrel{id}\to M\to 0$ is a free $\BF_2C_2$-resolution of $M$). To sum up, $\BF_2$-generators for the corner of $E_2$ are as follows:

\[
\xymatrix@=0pt{
{\phantom{4}} & x^3y+xy^3,x^4+y^4,x^2y^2   &x^2y^2z&  x^2y^2z^2  &x^2y^2z^3&  x^2y^2z^4   &x^2y^2z^5& x^2y^2z^6 \\
{\phantom{3}} & x^3+y^3,x^2y+xy^2   &0& 0    &0&0       &0& 0\\
{\phantom{2}} & x^2+y^2,xy   &xyz&  xyz^2  &xyz^3&  xyz^4   &xyz^5& xyz^6 \\
{\phantom{1}} & x+y   &0& 0    &0&0       &0& 0\\
{\phantom{0}} & 1   &z&  z^2  &z^3&  z^4   &z^5& z^6 \\
{\phantom{0}}\ar@{-}[]+U;[rrrrrrr]+UR\ar@{-}[]+R;[uuuuu]+UR  &  {\phantom{0}}  & {\phantom{1}} & {\phantom{2}}   & {\phantom{3}}&{\phantom{4}} &{\phantom{5}}&{\phantom{6}} &{\phantom{7}}
}
\]
A similar description of $E_2$ and the general fact that the spectral sequence of a wreath product collapses in $E_2$ may be found in \cite[IV, Theorem 1.7, p. 122]{Adem-Milgram}. So we have $E_\infty=E_2=\BF_2[z,\sigma_1,\sigma_2]/(z\sigma_1)$ where $\sigma_1=x+y$ and $\sigma_2=xy$ are the elementary symmetric polynomials. So according to Theorem \ref{thm:Carlsonliftinggensandrels}, we have $H^*(D_8;\BF_2)=\BF_2[w,\tau_1,\tau_2]/(R)$, where $|w|=1$, $|\tau_1|=1$, $|\tau_2|=2$ and $R$ is the lift of the relation $z\sigma_1=0$.

Recall that $|z|=(1,0)$ and that $|\sigma_1|=(0,1)$. So $z\sigma_1=0\in E_\infty^{1,1}=F^1H^2/F^2H^2$, where $F^*$ is some unknown filtration of $H^*=H^*(D_8;\BF_2)$. So $w\tau_1\in F^2H^2$ and we need to know what is $F^2H^2$. For total degree $2$, we have extension problems \eqref{equ:filtrationcohoonetotaldimension}:
\begin{align*}
F^1H^2\to &H^2 \to E_\infty^{0,2}=\langle x^2+y^2,xy\rangle=\BF_2\oplus \BF_2\\
F^2H^2\to &F^1H^2 \to E_\infty^{1,1}=0\\
0\to &F^2H^2 \to E_\infty^{2,0}=\langle z^2\rangle =\BF_2.
\end{align*}
So we deduce that $F^2H^2=F^1H^2=\langle z^2 \rangle\subset H^2(D_8;\BF_2)$. So we may take $w=z$ (see Remark \ref{rmk:horizontalaxisissubalgebraLHS}) and then $w\tau_1=\lambda w^2$ for some $\lambda\in \BF_2$. We cannot deduce the value of $\lambda$ from the spectral sequence. A computation with the bar resolution shows that $\lambda=0$, and hence $H^2(D_8;\BF_2)=\BF_2[w,\tau_1,\tau_2]/(w\tau_1)$.

\begin{exe}\label{exe:Poinacreseriesofdihedral}
Calculate the Poincar\'e series of $H^*(D_8;\BF_2)$, 
\[
P(t)=\sum_{i=0}^\infty \dim H^i(D_8;\BF_2)t^i,
\]
from the bigraded module $E_\infty$ of Example  \ref{example:dihedralaswreath}. Solution is in \ref{exe:Poinacreseriesofdihedralsolution}.
\end{exe}

\begin{exe}\label{exe:dihedralasC4:C2}
Fully describe the Lyndon-Hochschild-Serre spectral sequence with coefficients $\BF_2$ \ref{ss:LHScohomology} of the dihedral group $D_8$ described as a  semi-direct product:
\[
D_8=C_4\rtimes C_2=\langle a\rangle \rtimes \langle b\rangle,
\]
with ${}^ba=a^{-1}$. Recall that by Example \ref{example:dihedralaswreath} $H^*(D_8;\BF_2)=\BF_2[w,\tau_1,\tau_2]/(w\tau_1)$. Solution is in \ref{exe:dihedralasC4:C2solution}.
\end{exe}

\section{Edge morphisms}
\label{ss7:edgemorphisms}

Consider a cohomological type spectral sequence $E_2\Rightarrow H$ which is first quadrant, i.e., such that $E_2^{p,q}=0$ whenever $p<0$ or $q<0$. It is clear that there are monomorphisms and epimorphisms for all $p,q\geq 0$:
\begin{align*}
E_\infty^{0,q}\hookrightarrow\ldots \hookrightarrow E_r^{0,q}\hookrightarrow E_{r-1}^{0,q}\hookrightarrow\ldots\hookrightarrow E_3^{0,q}\hookrightarrow E_2^{0,q},\\
E_2^{p,0}\twoheadrightarrow E_3^{p,0}\twoheadrightarrow \ldots \twoheadrightarrow E_{r-1}^{p,0}\twoheadrightarrow E_r^{p,0}\twoheadrightarrow \ldots \twoheadrightarrow E_{\infty}^{p,0}.
\end{align*}
So $E_\infty^{0,q}$ is a $R$-submodule of $E_2^{0,q}$ and $E_\infty^{p,0}$ is a quotient $R$-module of $E_2^{p,0}$. Now, if we rewrite the extensions \eqref{equ:filtrationcohoonetotaldimension} together with $E_\infty^{a,b}\cong F^aH^{a+b}/F^{a+1}H^{a+b}$ we get
\begin{align}
&F^1H^{a+b}\to H^{a+b} \to E_\infty^{0,a+b}\nonumber\label{equ:filtrationcohoonetotaldimensionpq}\\
&F^2H^{a+b}\to F^1H^{a+b} \to E_\infty^{1,a+b-1}\nonumber\\
&\ldots\nonumber \\
&F^{a+b}H^{a+b}\to F^{a+b-1}H^{a+b}\to E_\infty^{a+b-1,1}\nonumber \\
&0\to F^{a+b}H^{a+b}\to E_\infty^{a+b,0}.\nonumber
\end{align}

We deduce that $E_\infty^{0,q}$ is a quotient module of $H^q=F^0H^q$ and that $F^pH^p=E_\infty^{p,0}$ is a submodule of $H^p$. Summing up, we have morphisms, called \emph{edge morphisms}:
\begin{align*}
H^q\twoheadrightarrow E_\infty^{0,q}\hookrightarrow E_2^{0,q},\\
E_2^{p,0}\twoheadrightarrow E_\infty^{p,0}\hookrightarrow H^p.
\end{align*}
For several spectral sequences, the edge morphisms may be explicitly described.  

For the Lyndon-Hochschild-Serre spectral sequence \ref{ss:LHScohomology} of a short exact sequence of groups, $N\stackrel{\iota}\to G\stackrel{\pi}\to Q$, we have $E_2^{0,q}=H^0(Q;H^q(N;R))=H^q(N;R)^Q$ and the edge morphism, \[
H^q(G;R)\twoheadrightarrow E_\infty^{0,q}\hookrightarrow E_2^{0,q}=H^q(N;R)^Q,
\] 
coincides with the restriction in cohomology 
\[
H^q(\iota)\colon H^q(G;R)\to H^q(N;R)^Q.
\]
We also have $E_2^{p,0}=H^p(Q;H^0(N;R))=H^p(Q;R^N)$ and the edge morphism,
\[
H^p(Q;R^N)=E_2^{p,0}\twoheadrightarrow E_\infty^{p,0}\hookrightarrow H^p(G;R),
\]
is exactly the inflation
\[
H^p(\pi)\colon H^p(Q;R^N)\to H^p(G;R).
\]

\begin{rmk}\label{rmk:horizontalaxisissubalgebraLHS}
We have seen that $E^{p,0}_\infty \subseteq H^p(G;R)$ as an $R$-module for all $p\geq 0$ but in fact $E^{*,0}_\infty\subseteq H^*(G;R)$ as a subalgebra because $E^{p,0}_\infty\cdot E^{p',0}_\infty\subseteq E^{p+p',0}_\infty$.
\end{rmk}

\begin{exe}\label{exe:C2C4C2F2finishoff}
Prove that in Example \ref{example:C2C4C2F2} we can deduce the algebra structure of $H^*(C_4;\BF_2)$ from the $E_\infty$-page. Solution in \ref{exe:C2C4C2F2finishoffsolution}. 
\end{exe}

Consider the Serre spectral sequence of a fibration $F\stackrel{\iota}\to E\stackrel{\pi}\to B$ \ref{ss:Serrecohomology} with $B$ simply connected and $F$ connected. The edge morphism 
\[
H^q(E;R)\twoheadrightarrow E_\infty^{0,q}\hookrightarrow E_2^{0,q}=H^q(F;R),
\] 
is
\[
H^q(\iota)\colon H^q(E;R)\to H^q(F;R).
\] 
The edge morphism
\[
H^p(B;R)=E_2^{p,0}\twoheadrightarrow E_\infty^{p,0}\hookrightarrow H^p(E;R),
\]
coincides with 
\[
H^p(\pi)\colon H^p(B;R)\to H^p(E;R).
\]
\begin{rmk}\label{rmk:horizontalaxeissubalgebraSerre}
We have seen that $E^{p,0}_\infty \subseteq H^p(E;R)$ as an $R$-module for all $p\geq 0$ but in fact $E^{*,0}_\infty\subseteq H^*(E;R)$ as a subalgebra because $E^{0,q}_\infty\cdot E^{0,q'}_\infty\subseteq E^{0,q+q'}_\infty$.
\end{rmk}

\begin{exe}\label{exe:edgemorphismshomologicalcase}
Figure out what are the edge morphisms for the Serre and the Lyndon-Hochschild-Serre homological spectral sequences.
\end{exe}

The next two examples show two particular situations where the edge morphisms give much information.

\begin{exm}[{\cite[p. 246]{Davis-Kirk}}]\label{example:orientedbordism}
We consider oriented bordism $\Omega_*^{SO}$, which is a cohomology theory. Then we have the Atiyah-Hirzebruch spectral sequence
\ref{ss:AtiyaHirzebruchcohomology} for the fibration  
\[
*\stackrel{\iota} \to X\stackrel{id} \to X,
\]
where $X$ is any topological space. It states that
\[
H_p(X;\Omega_q^{SO}(*))\Rightarrow \Omega_{p+q}^{´SO}(X).
\]
Also we have that 
\[
\Omega_q^{SO}(*)=\begin{cases}0\text{, $q<0$ or $q=1,2,3$,}\\ \BZ\text{, for $q=0,4$.}\end{cases}
\]
In this case, the edge morphism,
\[
\Omega_q^{SO}(*)\twoheadrightarrow H_0(X;\Omega_q^{SO}(*))=E^2_{0,q}\twoheadrightarrow E^\infty_{0,q} \hookrightarrow \Omega_q^{SO}(X),
\]
is given by $\Omega_q^{SO}(\iota)$. Note that the constant map $c\colon X\to *$ satisfies $c\circ \iota=1_*$. Hence, $\Omega_q^{SO}(*)\circ \Omega_q^{SO}(\iota)=1_{\Omega_q^{SO}(*)}$ and the edge morphism $\Omega_q^{SO}(\iota)$ is injective. In particular, $\Omega_q^{SO}(*)=H_0(X;\Omega_q^{SO}(*))$ and $E^2_{0,q}=E^{\infty}_{0,q}$ for all $q$. This last condition implies that all differentials arriving or emanating from the vertical axis must be zero:
\[
\xymatrix@=0pt{
4 & \Omega_4^{SO}(*) &H_1(X;\BZ)&H_2(X;\BZ)  &H_3(X;\BZ)& H_4(X;\BZ)\\
3 & 0   &0& 0    &0&0       \\
2 & 0   &0& 0    &0&0       \\
1 & 0   &0& 0    &0&0       \\
0 & \Omega_0^{SO}(*)   &H_1(X;\BZ)&H_2(X;\BZ)  &H_3(X;\BZ)& H_4(X;\BZ)\\
{\phantom{0}}\ar@{-}[]+U;[rrrrrr]+UR\ar@{-}[]+R;[uuuuu]+UR  &  0&1&2&3&4&{\phantom{5}}
}
\]
In particular, we deduce the oriented bordism of $X$ in low dimensions:
\[
\Omega_q^{SO}(X)=\begin{cases} H_q(X;\BZ)\text{, for $q=0,1,2,3$,}\\ \BZ\oplus H_4(X;\BZ)\text{, for $q=0,4$.}\end{cases}
\]
\end{exm}

\begin{exm}[{\cite[Proposition 7.3.2]{Evens}}]
\label{example:splitevens}
Consider an extension of groups $N\stackrel{\iota}\to G\stackrel{\pi}\to Q$ which is split, i.e., such that $G$ is the semi-direct product $G=N\rtimes Q$. We consider the Lyndon-Hochschild-Serre spectral sequence \ref{ss:LHScohomology} with coefficients in a ring $R$ with trivial $G$-action:
\[
E_2^{p,q}=H^p(Q;H^q(N;R))\Rightarrow H^{p+q}(G;R).
\]
We know that the edge morphism
\[
H^p(Q;R)=E_2^{p,0}\twoheadrightarrow E_\infty^{p,0}\hookrightarrow H^p(G;R),
\]
is exactly the inflation
\[
H^p(\pi)\colon H^p(Q;R)\to H^p(G;R).
\]
Because the extension is split, there is a homomorphism $s\colon Q\to G$ with $\pi\circ s=1_Q$. Hence, $H^p(s)\circ H^p(\pi)=1_{H^p(Q)}$, the inflation $H^p(\pi)$ is injective and we deduce that $E_2^{p,0}=E_\infty^{p,0}$. So all differentials arriving to the horizontal axis must be zero and we have an inclusion of algebras $H^*(Q;R)\subseteq H^*(G;R)$, see Remark \ref{rmk:horizontalaxisissubalgebraLHS}. This situation already occurred in Example \ref{example:dihedralaswreath} and appears again in Exercise \ref{exe:dihedralasC4:C2}.
\end{exm}

\section{Different spectral sequences with same target}
\label{ss8:diffsssametarget}

In Examples \ref{example:dihedralaswreath} and Exercise \ref{exe:dihedralasC4:C2} we have seen two different spectral sequences converging to the same target $H^*(D_8;\BF_2)$. In this section, we present another two examples of this phenomenon, the first with target $H^*(3^{1+2}_+;\BF_3)$ and the second with target $H_*(S^3;\BZ)$. From these examples, one sees that the price one pays for having an easily described $E_2$-page is more complicated differentials.

\begin{exm}\label{example:31+2+centralnoncentral}
We denote by $S=3^{1+2}_+$ the extraspecial group of order $27$ and exponent $3$. It has the following presentation
$$
S= \langle A,B,C | A^3=B^3=C^3=[A,C]=[B,C]=1\textit{, }[A,B]=C\rangle
$$
and it fits in the central extension
\[
\langle C\rangle =C_3 \rightarrow 3_+^{1+2} \rightarrow C_3\times C_3=\langle \bar A, \bar B\rangle.
\]
Leary describes in \cite{L1993} the Lyndon-Hochschild-Serre spectral sequence of this extension. Its second page is given by
\begin{equation}\label{equ:E2page31+2+central}
E^{*,*}_2=H^*(C_3;\BF_3)\otimes H^*(C_3\times C_3;\BF_3)=\Lambda(u)\otimes \BF_3[t]\otimes \Lambda(y_1,y
_2)\otimes \BF_3[x_1,x_2],
\end{equation}
with $|u|=|y_1|=|y_2|=1$ and $|t|=|x_1|=|x_2|=2$. The differentials in $E_*$ are the following, and the spectral sequence collapses at $E_6$.
\begin{enumerate}[(i)]
\item $d_2(u)=y_1y_2$, $d_2(t)=0$,
\item $d_3(t)=x_1y_2-x_2y_1$,
\item $d_4(t^iu(x_1y_2-x_2y_1))=it^{i-1}(x_1x_2^2y_2-x_1^2x_2y_1)$, $d_4(t^2y_i)=u(x_1y_2-x_2y_1)x_i$,\item $d_5(t^2(x_1y_2-x_2y_1))=x_1^3x_2-x_1x_2^3$,  $d_5(ut^2y_1y_2)= k u(x_1^3y_2-x_2^3y_1)$, $k \ne 0$.
\end{enumerate}
A long and intricate computation leads from $E_2$ to $E_6$ \cite{Diaz-Garaialde}. The following table contains representatives of classes that form an $\BF_3$-basis of $E_6^{n,m}$ for $0\leq n\leq 6$ and $0\leq m\leq 5$:  
{\small
$$
\xymatrix@=0pt{
&&&&&&&\\
5& & & & & & &\\
4& & & & & & &\\
3& & &uty_1y_2 &  & & &\\
2& &ty_1, ty_2 & &ty_1x_1, ty_1x_2&  &tx_1^2y_1, tx_1^2y_2&\\
& && &ty_2x_2&  &tx_2^2y_1, tx_2^2y_2&\\
1& &uy_1, uy_2 &uy_1y_2 &uy_1x_1, uy_1x_2 & &ux_1^2y_1, ux_1^2y_2&\\
& && &uy_2x_1, uy_2x_2  & &ux_2^2y_1, ux_2^2y_2 &\\
0&1&y_1,y_2&x_1,x_2 &y_1x_1, y_1x_2&x_1^2, x_2^2&x_1^2y_1, x_1^2y_2 &x_1^3, x_2^3\\
&&&  & y_2x_2 &x_1x_2 &x_2^2y_1, x_2^2y_2 &x_1^2x_2, x_1x_2^2 \\
&0 & 1 & 2 & 3 & 4 & 5 & 6 
\ar@{-}"1,2"+<-5pt,-2pt>;"10,2"+<-5pt,-5pt>;
\ar@{-}"1,2"+<5pt,-2pt>;"10,2"+<5pt,-5pt>;
\ar@{-}"1,3"+<17pt,-2pt>;"10,3"+<17pt,-5pt>;
\ar@{-}"1,4"+<15pt,-2pt>;"10,4"+<15pt,-5pt>;
\ar@{-}"1,5"+<28pt,-2pt>;"10,5"+<28pt,-5pt>;
\ar@{-}"1,6"+<15pt,-2pt>;"10,6"+<15pt,-5pt>;
\ar@{-}"1,7"+<27pt,-2pt>;"10,7"+<27pt,-5pt>;
\ar@{-}"1,8"+<25pt,-2pt>;"10,8"+<25pt,-5pt>;
\ar@{-}"1,2"+<-5pt,-2pt>;"1,8"+<25pt,-2pt>;
\ar@{-}"2,2"+<-5pt,-5pt>;"2,8"+<25pt,-5pt>;
\ar@{-}"3,2"+<-5pt,-5pt>;"3,8"+<25pt,-5pt>;
\ar@{-}"4,2"+<-5pt,-5pt>;"4,8"+<25pt,-5pt>;
\ar@{-}"6,2"+<-5pt,-5pt>;"6,8"+<25pt,-5pt>;
\ar@{-}"8,2"+<-5pt,-5pt>;"8,8"+<25pt,-5pt>;
\ar@{-}"10,2"+<-5pt,-5pt>;"10,8"+<25pt,-5pt>;
}
$$
}

It turns out that this description of the corner of $E_6$ determines the rest of $E_6$ as there are both vertical and horizontal periodicities. More precisely, there are $\BF_3$-isomorphisms  $E_6^{n,m}\cong E_6^{n,m+6}$ for $n,m\geq 0$ and $E_6^{n,m}\cong E_6^{n+2,m}$ for $n\geq 5$ and $m\geq 0$. 
The extraspecial group  $S=3^{1+2}_+$ also fits in an extension
\[
\langle B,C\rangle =C_3\times C_3 \rightarrow 3_+^{1+2} \rightarrow C_3=\langle \bar A\rangle,
\]
where the action is given by the matrix $\left(\begin{smallmatrix} 1 &0\\ 1& 1 \end{smallmatrix}\right)$. In this case, the page $E_2^{*,*}=H^*(C_3;H^*(C_3\times C_3;\BF_3))$ is quite complicated as the action of $C_3$ is far from being trivial. It is described in \cite{Siegel}, where the author also shows that $E_2=E_\infty$. The page $E_2$ has $9$ generators, $x_1,\gamma_1,x_2,y_2,\gamma_2,x_3,x_6,z_2,z_3$, and they lie in the following positions:
\[
\xymatrix@=0pt{
6  & x_6  &   &    &  \\
5  &    &    &    &  \\
4  &    &    &    &  \\
3  & x_3   &    &    &  \\
2  & x_2,y_2   &  z_3  &     &  \\
1  & x_1   &  z_2  &     &  \\
0  &    &  \gamma_1  &  \gamma_2   &  \\
{\phantom{0}}\ar@{-}[]+U;[rrrr]+UR\ar@{-}[]+R;[uuuuuuu]+UR  & 0&1&2&3&
}
\]Compare to the $E_2$-page \eqref{equ:E2page31+2+central}, which has $6$ generators all of which lie on the axes.
\end{exm}

\begin{exm}\label{example:C_2S_3RP3}
In Example \ref{example:HopffibrationS1S3S2}, we saw that the Hopf fibration $S^1\to S^3\to S^2$ gives rise to a spectral sequence converging to $H_*(S^3)$. Consider now the fibration
\[
C_2\to S^3\to \BR\BP^3,
\]
where $C_2$ acts by the antipodal map and $\BR\BP^3$ is the projective space of dimension $3$. Although the fiber $C_2$ is not connected, we have the associated Serre spectral sequence \ref{ss:Serrehomology}:
\[
E^2_{p,q}=H_p(\BR\BP^3;H_q(C_2;\BZ))\Rightarrow H_{p+q}(S^3;\BZ).
\]
Of course now we have $H_0(C_2;\BZ)=\BZ\oplus \BZ$ with $C_2=\pi_1(\BR P^3)$ interchanging the two $\BZ$'s, and $H_q(C_2;\BZ)=0$ for $q>0$. Hence, $E^2_{p,0}$ is the twisted homology $H_p(\BR\BP^3;\BZ\oplus \BZ)$ and $E^2_{p,q}=0$ for $q>0$. Then the spectral sequence is concentrated at the horizontal axis and collapses at the $E^2$-page. We deduce that the twisted homology $H_p(\BR\BP^3;\BZ\oplus \BZ)$ must be $\BZ$ for $p=0,3$ and $0$ otherwise. 

Let us directly compute the twisted homology $H_p(\BR\BP^3;\BZ\oplus \BZ)$ via the universal cover $S^3$ of $\BR\BP^3$ \cite[3.H]{Hatcher}. The sphere $S^3$ has a structure of free $C_2$-complex with cells 
\[
S^3=e^3_+\cup e^3_-\cup e^2_+\cup e^2_-\cup e^1_+\cup e^1_-\cup e^0_+\cup e^0_-,
\]
where $C_2$ interchanges each pair of cells on every dimension. This gives rise to the following  chain complex over $\BZ C_2$:
\[
0\to \BZ\oplus \BZ\to \BZ\oplus \BZ\to \BZ\oplus \BZ\to \BZ\oplus \BZ\to \BZ\to 0,
\]
with $d(1,0)=(1,-1)$, $d(0,1)=(-1,1)$ in dimensions $1$ and $3$, $d(1,0)=d(0,1)=(1,1)$ in dimension $2$ and $d(1,0)=d(0,1)=1$ in dimension $0$. Set $M$ to be equal to the $\BZ C_2$-module $H_0(C_2;\BZ)$. Then we need to compute the homology of
\[
0\to (\BZ\oplus \BZ)\otimes_{C_2} M \to (\BZ\oplus \BZ)\otimes_{C_2} M \to (\BZ\oplus \BZ)\otimes_{C_2} M \to (\BZ\oplus \BZ)\otimes_{C_2} M \to 0.
\]
Setting $a=(1,0)\otimes(1,0)=(0,1)\otimes(0,1)$ and $b=(1,0)\otimes(0,1)=(0,1)\otimes(1,0)$ on every dimension, it turns out that $d(a)=a-b$, $d(b)=b-a$ in dimensions $1$ and $3$, $d(a)=d(b)=a+b$ in dimension $2$ and $d(a)=d(b)=0$ in dimension $0$. Thus we obtain the desired homology $H_*(S^3;\BZ)$.
\end{exm}

\begin{exe}\label{exe:IcoS_3PoincareS^3}
Do similar computations to those of Example \ref{example:C_2S_3RP3} for the fibration
\[
I\to S^3 \to PS^3,
\]
where $I$ is the binary icosahedral group of order $120$ and $PS^3$ is the Poincar\'e homology sphere. This space satisfies $H_p(PS^3;\BZ)=H_p(S^3;\BZ)$ for all $p$ and $\pi_1(PS^3)=I$, which is thus a perfect group.
\end{exe}

\section{A glimpse into the black box}
\label{ss9:glimpse}
In this section, we show how a filtration gives rise to a spectral sequence. Spectral sequences may also be constructed from exact couples, and both approaches are equivalent.  So let $(C,d)$ be a chain complex and let $F_*C$ be an increasing filtration of $(C,d)$, i.e., an ordered family of chain subcomplexes:
\[
\ldots \subseteq F_{n-1}C\subseteq F_nC\subseteq \ldots \subseteq C.
\]
We may define a filtration of the homology $H$ of $C$ via the inclusions: $F_nH=\Img(H(F_nC,d)\longrightarrow H(C,d))$. Thus we have an ordered family of graded modules:
\[
\ldots \subseteq F_{n-1}H\subseteq F_nH\subseteq \ldots \subseteq H.
\]
On each total degree $p+q$, we obtain a filtration of the modules $C_{p+q}$ and $H_{p+q}$:
\begin{align*}
\ldots &\subseteq F_{p-1}C_{p+q}\subseteq F_pC_{p+q}\subseteq \ldots \subseteq C_{p+q}\\
\ldots &\subseteq F_{p-1}H_{p+q}\subseteq F_pH_{p+q}\subseteq \ldots \subseteq H_{p+q}.
\end{align*}
We will assume that the filtration is bounded, i.e., that for each total degree $n$, there exist $s$ and $r$ such that $F_sC_n=0$ and $F_tC_n=C_n$. Then we have:
\begin{align}
0=F_sC_{p+q}&\subseteq \ldots\subseteq F_{p-1}C_{p+q}\subseteq F_pC_{p+q}\subseteq \ldots \subseteq F_rC_{p+q}=C_{p+q}\label{equ:constructionfiltC}\\
0=F_sH_{p+q}&\subseteq \ldots\subseteq F_{p-1}H_{p+q}\subseteq F_pH_{p+q}\subseteq \ldots \subseteq F_rH_{p+q}=H_{p+q}.\label{equ:constructionfiltH}
\end{align}
Because of the boundedness assumption, there exists  a spectral sequence $\{E^r_{*,*},d_r\}_{r\geq 0}$ that converges to $H$ in the sense explained in \eqref{enumerate:convergencestabilizes} and \eqref{enumerate:convergencefiniteextensions}. More precisely, $E^0_{*,*}$ is obtained by taking quotients on \eqref{equ:constructionfiltC}:
\[
E^0_{p,q}=F_pC_{p+q}/F_{p-1}C_{p+q},
\]
and the stable term $E^\infty_{p,q}$ at $(p,q)$ is isomorphic to considering quotients on \eqref{equ:constructionfiltH}:
\[
E^\infty_{p,q}\cong F_pH_{p+q}/F_{p-1}H_{p+q}.
\]
This correspond to the Equations \eqref{equ:filtrationhoonetotaldimension}. So the problem is to define the rest of pages $E^r_{*,*}$ and their differentials and see how they fit together. 

Define $Z^r_{p,q}=F_pC_{p+q}\cap d^{-1}(F_{p-r}C_{p+q-1})$ and $B^r_{p,q}=F_pC_{p+q}\cap d(F_{p+r}C_{p+q+1})$:
\[
\xymatrix@R=0pt@C=10pt{
&&&          &&&F_{p+r}C_{p+q+1}\ar[llld]_d\\   
&&&F_pC_{p+q}\ar[llld]_d&&&\\
F_{p-r}C_{p+q-1}\\
}
\]
Then we set
\[
E^r_{p,q}=\frac{Z^r_{p,q}}{Z^{r-1}_{p-1,q+1}+B^{r-1}_{p,q}},
\]
and one has to check that:
\begin{enumerate}[(a)]
\setcounter{enumi}{4}
\item The differential $d$ induces a differential $d^r_{p,q}\colon E^r_{p,q}\to E^r_{p-r,q+r-1}$.
\item When taking homology we have $H_{p,q}(E^r_{*,*},d^r)\cong E^{r+1}_{p,q}$.
\end{enumerate}
These two facts are consequence of intricate diagram chasings, and it is here where things become  delicate. We do not reproduce the details (that may be found in \cite[2.2]{McCleary} or \cite[5.4]{Weibel} for instance), but instead we explain the crucial step. 

The aforementioned chasing arguments reduce to consider the following situation: We are given a $3\times 3$-diagram of modules and inclusions,
\[
\xymatrix@=10pt{A\ar@{^(->}[r]\ar[d]^d&B\ar[d]^d\ar@{^(->}[r]&C\ar[d]^d\\
D\ar[d]^d\ar@{^(->}[r]&E\ar[d]^d\ar@{^(->}[r]&F\ar[d]^d\\
G\ar@{^(->}[r]&H\ar@{^(->}[r]&I,
}
\]
where $d\circ d=0$, and we need to prove \eqref{equ:3times3one} and \eqref{equ:3times3two} below. We define $E_1$ as taking quotient and homology (in this order):
\begin{align*}
E_1^E&=H(B/A\to E/D\to H/G),\\
E_1^G&=H(D\to G\to 0),\\
E_1^C&=H(0\to C/B\to F/E),
\end{align*}
and $E_\infty$ as taking homology and quotient (in this order):
\[
E^E_\infty=\frac{\Img(H(B\to E\to H)\to H(C\to F\to I))}{\Img(H(A\to D\to G)\to H(C\to F\to I))}.
\]
Then the following hold:
\begin{enumerate}[(a)]
\setcounter{enumi}{6}
\item The differential $d$ induces a differential $E_1^C\stackrel{d^1}\to E_1^E\stackrel{d^1}\to E_1^G$.\label{equ:3times3one}
\item When taking homology we have $H(E_1^C\stackrel{d^1}\to E_1^E\stackrel{d^1}\to E_1^G)\cong E^E_\infty$.\label{equ:3times3two}
\end{enumerate}

\begin{exe}
Prove that \eqref{equ:3times3one} and \eqref{equ:3times3two} hold.
\end{exe}

The situation for a decreasing filtration of a cochain complex is similar. Considering a product in the cochain complex gives rise to a spectral sequence of algebras, see \cite[2.3]{McCleary}.

\section{Appendix I}
\label{section:AppendixI}

\begin{thm}[Serre spectral sequence for homology, {\cite[Theorem 5.1]{McCleary}}]
\label{ss:Serrehomology}
Let $M$ be an abelian group and let $F\to E\to B$ be a fibration with $B$ path-connected. Then there is a first quadrant homological spectral sequence
\[
E^2_{p,q}=H_p(B;H_q(F;M))\Rightarrow H_{p+q}(E;M).
\]
\end{thm}

\begin{thm}[Serre spectral sequence for cohomology, {\cite[Theorem 5.2]{McCleary}}]
\label{ss:Serrecohomology}
Let $R$ be a ring and let $F\to E\to B$ be a fibration with $B$ path-connected. Then there is a first quadrant cohomological spectral sequence of algebras and converging as an algebra
\[
E_2^{p,q}=H^p(B;H^q(F;R))\Rightarrow H^{p+q}(E;R).
\]
\end{thm}

In the latter theorem, if $\pi_1(B)=0$ and $R$ is a field we have
\[
E_2^{p,q}=H^p(B;R)\otimes H^q(F;R))\Rightarrow H^{p+q}(E;R).
\]

\begin{proof}[Sketch of proof of Theorem \ref{ss:Serrecohomology}]
We follow \cite[Chapter 5]{McCleary}. Denote by $F\to E\stackrel{\pi}\to B$ the given fibration and consider the skeletal filtration of $B$
\[
\emptyset \subseteq B^0\subseteq B^1\subseteq B^2 \ldots \subseteq B.
\]
Then we have a filtration of $E$ given by $J^s=\pi^{-1}(B^s)$
\[
\emptyset \subseteq J^0\subseteq J^1\subseteq J^2 \ldots \subseteq E.
\]
Consider the cochain complex of singular cochains on $E$ with coefficients in $R$, $C^*(E;R)$, given by
\[
C^s(E;R)=\{\text{functions $f\colon C_s(E)\to R$}\},
\]
where $C_s(E)=\{\sigma\colon \Delta^s\to E\text{ continuous}\}$ are the singular $s$-chains in $E$. Now define a decreasing filtration of $C^*(E;R)$ by
\[
F^sC^*(E;R)=\ker(C^*(E;R)\to C^*(J^{s-1};R)),
\]
i.e., the singular cochains of $E$ that vanish on chains in $J^{s-1}$. It turns out that 
\[
E_1^{p,q}=H^{p+q}(J^p,J^{p-1};R),
\]
the relative cohomology of the pair $(J^p,J^{p-1})$. 

\begin{center}
\begin{tikzpicture}[scale=0.5,line cap=round,line join=round,>=triangle 45,x=1.0cm,y=1.0cm]
\clip(4.,4.) rectangle (19.,14.);
\draw [rotate around={0.:(10.196528644356324,8.745706182643223)}] (10.196528644356324,8.745706182643223) ellipse (4.023258888640295cm and 0.5598356908667658cm);
\draw [rotate around={0.:(10.423414634146365,8.800158820192856)}] (10.423414634146365,8.800158820192856) ellipse (3.0946691589261404cm and 0.4091816868677692cm);
\draw [rotate around={-0.29841281999487734:(10.949790130459338,8.79108338060125)}] (10.949790130459338,8.79108338060125) ellipse (1.7628103981926149cm and 0.2667700256708816cm);
\draw [rotate around={0.:(10.287283040272275,13.301576857629039)}] (10.287283040272275,13.301576857629039) ellipse (4.0232588886403855cm and 0.5598356908667783cm);
\draw [rotate around={0.:(10.514169030062511,13.356029495178673)}] (10.514169030062511,13.356029495178673) ellipse (3.0946691589264694cm and 0.4091816868678127cm);
\draw [rotate around={-0.29841281999487734:(11.040544526375733,13.346954055587062)}] (11.040544526375733,13.346954055587062) ellipse (1.7628103981922751cm and 0.2667700256708302cm);
\draw (6.265040213902779,13.314157978399919)-- (6.178560443602286,8.774407368259995);
\draw (14.31054192891266,13.301576857629039)-- (14.21978753299662,8.745706182643223);
\draw (7.42158869716836,13.340998009355824)-- (7.328745475220225,8.800158820192856);
\draw (13.605096389090217,13.335913893144655)-- (13.514341993173801,8.820274422226719);
\draw (9.279733487629434,13.36854276243051)-- (9.18775657656571,8.792251296030193);
\draw (12.78498510346572,13.299688015536473)-- (12.712576619492141,8.78190220097087);
\draw [rotate around={0.6015083943170401:(10.239228587634733,5.8297901304594415)}] (10.239228587634733,5.8297901304594415) ellipse (3.816102698635533cm and 1.614106223405729cm);
\draw [rotate around={-0.9911775750107868:(10.622802041973896,5.954622802041972)}] (10.622802041973896,5.954622802041972) ellipse (2.8692441132799704cm and 1.1625025333799717cm);
\draw [rotate around={-0.6700955222245235:(10.926193987521263,5.960272263187744)}] (10.926193987521263,5.960272263187744) ellipse (1.7464385694758342cm and 0.8008271251531659cm);
\draw (17.17917186613727,13.456585365853655)-- (17.12471922858764,8.573998865570047);
\draw[color=black] (7,5.8) node {$B$};
\draw[color=black] (8.6,5.8) node {$B^s$};
\draw[color=black] (10.5,5.8) node {$B^{s-1}$};
\draw[color=black] (6.7,10.8) node {$E$};
\draw[color=black] (8.6,10.8) node {$J^s$};
\draw[color=black] (10.5,10.8) node {$J^{s-1}$};
\draw[color=black] (16.5,10.8) node {$F$};

\end{tikzpicture}
\end{center}

The keystep of the proof relies in showing that 
\[
H^{p+q}(J^p,J^{p-1};R)\cong C^p(B;H^q(F;R)),
\]
the twisted $p$-cochains of $B$ with coefficients in the $q$-cohomology of the fiber, which are defined as follows:
\[
\{\text{functions $f\colon C_p(B)\to \bigcup_{b\in B} H^q(\pi^{-1}(b);R)$ such that $f(\sigma)\in H^q(\pi^{-1}(\sigma(v_0));R)$}\}.
\]

Note that, as $B$ is connected, $\pi^{-1}(b)$ has the same homotopy type ($F$) for all $b\in B$. Also, we assume that $\Delta^p\subseteq \BR^p$ is spanned by the vertices $\{v_0,\ldots,v_{p+1}\}$. The (twisted) differential $C^p(B;H^q(F;R))\to C^{p+1}(B;H^q(F;R))$ is given by
\[
\partial(f)(\sigma)=H^q(\Phi;R)(f(\sigma_0))+\sum_{i=1}^{p+1}(-1)^i f(\sigma_i),
\]
where $\sigma_i$ is the restriction of $\sigma$ to the $i$-th face of $\Delta^p$, $\{v_0,\ldots,\hat {v_i},\ldots,v_{p+1}\}$, and the map $\Phi\colon \pi^{-1}(\sigma(v_0))\to \pi^{-1}(\sigma(v_1))$ is a ``lift'' of the path $\sigma_{[v_0,v_1]}$ starting in $v_0$ and ending in $v_1$.
\begin{center}
\definecolor{zzttqq}{rgb}{0.26666666666666666,0.26666666666666666,0.26666666666666666}
\begin{tikzpicture}[scale=0.5,line cap=round,line join=round,>=triangle 45,x=1.0cm,y=1.0cm]
\clip(0.,3.) rectangle (15.,14.);
\fill[color=zzttqq,fill=zzttqq,fill opacity=0.1] (8.,6.) -- (11.461644923425984,6.7952127056154215) -- (12.078774815655139,5.143482699943272) -- cycle;
\fill[color=zzttqq,fill=zzttqq,fill opacity=0.1] (8.085581395348841,13.583641520136123) -- (11.3164378899603,13.728848553601805) -- (11.26218765339866,8.864410675241322) -- (8.03112875779921,8.701055019852516) -- cycle;
\fill[color=zzttqq,fill=zzttqq,fill opacity=0.1] (3.3119001701644963,6.867816222348263) -- (1.4605104934770319,5.161633579126483) -- (3.929030062393651,4.471900170164487) -- cycle;
\draw [rotate around={0.:(10.196528644356324,8.745706182643223)}] (10.196528644356324,8.745706182643223) ellipse (4.023258888640295cm and 0.5598356908667658cm);
\draw [rotate around={0.:(10.287283040272275,13.301576857629039)}] (10.287283040272275,13.301576857629039) ellipse (4.0232588886403855cm and 0.5598356908667783cm);
\draw (6.265040213902779,13.314157978399919)-- (6.178560443602286,8.774407368259995);
\draw (14.31054192891266,13.301576857629039)-- (14.21978753299662,8.745706182643223);
\draw [rotate around={0.6015083943170401:(10.239228587634733,5.8297901304594415)}] (10.239228587634733,5.8297901304594415) ellipse (3.816102698635533cm and 1.614106223405729cm);
\draw [color=zzttqq] (8.,6.)-- (11.461644923425984,6.7952127056154215);
\draw [color=zzttqq] (11.461644923425984,6.7952127056154215)-- (12.078774815655139,5.143482699943272);
\draw [color=zzttqq] (12.078774815655139,5.143482699943272)-- (8.,6.);
\draw (8.085581395348841,13.583641520136123)-- (8.03112875779921,8.701055019852515);
\draw (11.3164378899603,13.728848553601805)-- (11.26198525241067,8.846262053318197);
\draw [color=zzttqq] (8.085581395348841,13.583641520136123)-- (11.3164378899603,13.728848553601805);
\draw [color=zzttqq] (11.3164378899603,13.728848553601805)-- (11.26218765339866,8.864410675241322);
\draw [color=zzttqq] (11.26218765339866,8.864410675241322)-- (8.03112875779921,8.701055019852516);
\draw [color=zzttqq] (8.03112875779921,8.701055019852516)-- (8.085581395348841,13.583641520136123);
\draw [color=zzttqq] (3.3119001701644963,6.867816222348263)-- (1.4605104934770319,5.161633579126483);
\draw [color=zzttqq] (1.4605104934770319,5.161633579126483)-- (3.929030062393651,4.471900170164487);
\draw [color=zzttqq] (3.929030062393651,4.471900170164487)-- (3.3119001701644963,6.867816222348263);
\draw[color=black] (3,5.8) node {$\Delta^2$};
\draw[color=black] (1,5.1) node {$v_0$};
\draw[color=black] (3.2,7.2) node {$v_1$};
\draw[color=black] (4,4) node {$v_2$};
\draw[color=black] (5,5.8) node {$\stackrel{\sigma}\longrightarrow$};
\draw[color=black] (7.3,6) node {$\sigma(v_0)$};
\draw[color=black] (11,7.2) node {$\sigma(v_1)$};
\draw[color=black] (12.5,4.8) node {$\sigma(v_2)$};
\draw[color=black] (6.5,10.8) node {$\pi^{-1}(\sigma(v_0))$};
\draw[color=black] (13,10.8) node {$\pi^{-1}(\sigma(v_1))$};
\draw[color=black] (9.5,10.8) node {$\stackrel{\Phi}\longrightarrow$};
\end{tikzpicture}
\end{center}
Then $E_2^{*,*}$ is obtained as the cohomology of $C^p(B;H^q(F;R))$, and this gives the description in the statement of the theorem.
\end{proof}
\begin{thm}[Lyndon-Hochschild-Serre spectral sequence for homology, {\cite[6.8.2]{Weibel}}]
\label{ss:LHShomology}
Let $N\to G\to Q$ be a short exact sequence of groups and let $M$ be a $G$-module. Then there is a first quadrant homological spectral sequence
\[
E^2_{p,q}=H_p(Q;H_q(N;M))\Rightarrow H_{p+q}(G;M).
\]
\end{thm}

\begin{thm}[Lyndon-Hochschild-Serre spectral sequence for cohomology, {\cite[6.8.2]{Weibel}}]
\label{ss:LHScohomology}
Let $N\to G\to Q$ be a short exact sequence of groups and let $R$ be a ring. Then there is a first quadrant cohomological spectral sequence of algebras converging as an algebra
\[
E_2^{p,q}=H^p(Q;H^q(N;R))\Rightarrow H^{p+q}(G;R).
\]
\end{thm}

In the latter theorem, if $N\leq Z(G)$ and $R$ is a field we have
\[
E_2^{p,q}=H^p(Q;R)\otimes H^q(N;R)\Rightarrow H^{p+q}(G;R).
\]

\begin{proof}[Sketch of proof of Theorem \ref{ss:LHScohomology}]
We follow \cite[XI, Theorem 10.1]{MacLane} but consult \cite{Hochschild-Serre} for a direct filtration on the bar resolution of $G$. We do not provide details about the multiplicative structure. Consider the double complex
\[
C^{p,q}=\Hom_{RQ}(B_p(Q),\Hom_{RN}(B_q(G),R)\cong \Hom_{RG}(B_p(Q)\otimes B_q(G),R),
\]
where $B_*(Q)$, $B_*(G)$ denote the corresponding bar resolutions and $RQ$, $RN$ and $RG$ denote group rings. So $B_p(Q)$ is the free $R$-module with basis $Q^{p+1}$ and with differential the $R$-linear extension of $\partial(q_0,\ldots,q_p)=\sum_{i=0}^p (-1)^i(q_0,\ldots,\hat q_i,\ldots,q_p)$. Then the horizontal and vertical differentials of $C^{*,*}$ are given by
\begin{align*}
&\partial_h(f)(b\otimes b')=(-1)^{p+q+1}f(\partial(b)\otimes b')\text{, and}\\
&\partial_v(f)(b\otimes b')=(-1)^{q+1}f(b\otimes \partial(b')).
\end{align*}
There are two filtrations of the graded differential $R$-module $(\total(C),\partial_h+\partial_v)$ obtained by considering either all rows above a given row or all columns to the right of a given column:

\begin{minipage}{0.5\textwidth}
\begin{align*}
&\xymatrix@=10pt{
 \circ & \bullet\ar@{-}[]+C;[rrrddd]+C & \circ & \circ & \circ \\
 \circ & \frm{**}\circ & \bullet \ar[]+C;[r]+C^{d_h}\ar[]+C;[u]+C^{d_v}& \circ & \circ \\
 \circ \ar@{--}[]+C;[rrrr]+C & \circ & \circ & \bullet & \circ \\
 \cdot & \cdot & \cdot & \cdot & \cdot \\
 \cdot & \cdot & \cdot & \cdot & \cdot \\
  }
&\\
&\hspace{10pt}\text{Filtration by rows}
\end{align*}
\end{minipage}
\begin{minipage}{0.5\textwidth}
\begin{align*}
&\xymatrix@=10pt{
 \cdot & \cdot\ar@{-}[]+C;[rrrddd]+C & \circ \ar@{--}[]+C;[dddd]+C& \circ & \circ \\
 \cdot & \cdot & \bullet & \circ & \circ \\
 \cdot & \cdot & \circ & \bullet\ar[]+C;[r]+C^{d_h}\ar[]+C;[u]+C^{d_v} & \circ \\
 \cdot & \cdot & \circ & \circ & \bullet \\
 \cdot & \cdot & \circ & \circ & \circ  \\}
&\\
&\text{Filtration by columns}
\end{align*}
\end{minipage}

\vspace{10pt}These two filtrations give rise to spectral sequences converging to the cohomology of $(\total(C),\partial_h+\partial_v)$ \cite[Section 5.6]{Weibel},
\begin{align*}
&{}^rE_2^{p,q}=H^p_vH^q_h(C)\Rightarrow H^{p+q}(\total(C))\text{, and}\\
&{}^cE_2^{p,q}=H^p_hH^q_v(C)\Rightarrow H^{p+q}(\total(C)),
\end{align*}
where the subscripts $h$ and $v$ denote taking cohomology with respect to $\partial_h$ or $\partial_v$ respectively. For the former spectral sequence, we obtain ${}^rE_2^{p,q}=H^p(G;R)$ for $q=0$ and $E_2^{p,q}=0$ for $q>0$ by \cite[XI, Lemma 9.3]{MacLane}. For the latter spectral sequence we get
\[
{}^cE_1^{p,q}=\Hom_{RQ}(B_p(Q),H^q(N;R)),
\]
as $B_p(Q)$ is a free $RQ$-module and hence commutes with cohomology, and then
\[
{}^cE_2^{p,q}=H^p(Q;H^q(N;R)).
\]
\end{proof}
\begin{thm}[Atiyah-Hirzebruch spectral sequence, {\cite[Theorem 9.22]{Davis-Kirk}}]
\label{ss:AtiyaHirzebruchcohomology}
Let $h$ be a cohomology theory and let $F\to E\to B$ be a fibration with $B$ path-connected. Assume $h^q(F)=0$ for $q$ small enough. Then there is a ``half-plane'' cohomological spectral sequence
\[
E_2^{p,q}=H^p(B;h^q(F))\Rightarrow h^{p+q}(E).
\]
\end{thm}

\section{Appendix II}
\label{section:AppendixII}

\subsection{Homotopy groups of spheres.} 
The Serre spectral sequence may be used to obtain general results about homotopy groups of spheres. For instance, the cohomology Serre spectral sequence with coefficients $\BQ$ is the tool needed for the following result.
\begin{thm}[{\cite[Theorem 1.21]{SSch1}\cite[Theorem 10.10]{Davis-Kirk}}]
The groups $\pi_i(S^n)$ are finite for $i > n$, except for $\pi_{4n-1}(S^{2n})$ which is the direct sum of $\BZ$ with a finite group.
\end{thm}

Using as coefficients the integers localized at $p$, $\BZ_{(p)}$, i.e., the subring of $\BQ$ consisting
of fractions with denominator relatively prime to $p$, yields the next result.
\begin{thm}[{\cite[Theorem 1.28]{SSch1}\cite[Corollary 10.13]{Davis-Kirk}}]
For $n\geq 3$ and $p$ a prime, the $p$-torsion subgroup of $\pi_i(S^n)$ is zero for $i < n+2p-3$ and $C_p$ for $i=n+2p-3$.
\end{thm}

The Serre spectral sequence may be further exploited to compute more homotopy group of spheres \cite[Chapter 12]{Mosher-Tangora} but it does not give a full answer. The EHP spectral sequence is another tool to compute homotopy group of spheres, see for instance \cite[p. 43]{SSch1}.

\subsection{Stable homotopy groups of spheres.}
\label{subsection:Adamsspectralsequence}
These groups are defined as follows:
\[
\pi^S_k=\varinjlim_{n\to \infty} (\ldots\to \pi_{n+k}(S^n)\to\pi_{n+k+1}(S^{n+1})\to\ldots),
\]
where the homomorphisms are given by suspension:
\[
S^{n+k}\stackrel{f}\longrightarrow S^n \leadsto \Sigma S^{n+k}=S^{n+k+1}\stackrel{\Sigma f}\longrightarrow S^{n+1}=\Sigma S^n.
\]
In fact, by Freudenthal's suspension theorem, $\pi^S_k=\pi_{n+k}(S^n)$ for $n>k+1$, i.e., all morphisms become isomorphism for $n$ large enough. The Adams spectral sequence is the tool to compute  stable homotopy groups of spheres: For each prime $p$, there is a spectral sequence which second page is 
\[
E_2^{s,t}=\Ext^{s,t}_{\CA_p}(\BF_p,\BF_p)
\]
and converging to $\pi^S_*$ modulo torsion of order prime to $p$. Here, $\CA_p$ is the Steenrod algebra. The $E_2$-page is so complicated that the  May spectral sequence is used to determine it. See \cite[Chapter 18]{Mosher-Tangora}, \cite{SSch2} and \cite[Chapter 9]{McCleary}.

\subsection{Cohomology operations and Steenrod algebra.}
Cohomology operations of type $(n,m,\BF_p,\BF_p)$ are exactly the  natural transformations between the functors $H^n(-;\BF_p)$ and $H^m(-;\BF_p)$:
\[
\theta\colon H^n(-;\BF_p)\Rightarrow H^m(-;\BF_p).
\]
Among cohomology operations, we find stable cohomology operations, i.e., those families of cohomology operations of fixed degree $m$, 
\[
\{\theta_n\colon H^n(-;\BF_p)\rightarrow H^{n+m}(-;\BF_p)\}_{n\geq 0},
\] 
that commute with the suspension isomorphism:
\[
\xymatrix{H^n(X;\BF_p)\ar[r]^{\theta_n}\ar[d]^\Sigma & H^{n+m}(X;\BF_p)\ar[d]^\Sigma\\
H^{n+1}(\Sigma X;\BF_p)\ar[r]^{\theta_{n+1}} & H^{n+m+1}(\Sigma X;\BF_p).}
\]
These stable cohomology operations for mod $p$ cohomology are assembled together to form the Steenrod algebra $\CA_p$. Moreover, cohomology operations of a given type $(n,m,\BF_p,\BF_p)$ are in bijection with the cohomology group $H^m(K(\BF_p,n);\BF_p)$. The cohomology ring $H^*(K(\BF_p,n);\BF_p)$ may be determined via the Serre spectral sequence for all $n$ and for all primes $p$, although this computation is much harder than the already seen in Example \ref{example:cohoK(Z,n)}. Hence, information about $H^*(K(\BF_p,n);\BF_p)$ provides insight into the Steenrod algebra $\CA_p$, which in turn is needed in the Adams spectral sequence \ref{subsection:Adamsspectralsequence}. See \cite[Theorem 1.32]{SSch1},
\cite[Theorem 6.19]{McCleary},
\cite[Chapter 9]{Mosher-Tangora} and
\cite[Section 10.5]{Davis-Kirk}.

\subsection{Hopf invariant one problem}
Adams invented and used his spectral sequence to solve this problem (although later on he gave another proof via secondary cohomology operations). The Hopf invariant of a map $S^{2n-1}\stackrel{f}\to S^n$ is the only integer $H(f)$ satisfying:
\[
x^2=H(f)y,
\]
where $x$ and $y$ are the generators in degrees $n$ and $2n$ of the integral cohomology $H^*(X;\BZ)$ of certain space $X$. This space obtained by attaching a $2n$-dimensional cell to $S^n$ via $f$:
\[
X=S^n\cup_f D^{2n}.
\]
The Hopf invariant one problem consists of determining for which values of $n$ there exists a map $f$ with $H(f)=1$, and its solution is that  $n\in\{1,2,4,8\}$. The problem can also be phrased as determining for which $n$ does $\BR^n$ admit a division algebra structure. 
%This is equivalent to that the only $n$-spheres which are $H$-spaces are $S^0$, $S^1$, $S^3$, and $S^7$ and also to that the only spherical fibrations $S^l\to S^m\to S^n$ occur for $(l,m,n)\in \{(0,1,1),(1,3,2),(3,7,4),(7,15,8)\}$ \cite[4.B]{Hatcher}. 
See \cite[4.B]{Hatcher} and \cite[Theorem 9.38]{McCleary} or Adam's original papers \cite{Adams58}, \cite{Adams60}. 

\subsection{Segal's conjecture.}
This conjecture states that for any finite group $G$, the natural map from the completion of the Burnside ring of $G$ at its augmentation ideal to the stable cohomotopy of its classifying space is an isomorphism:
\[
{A(G)^{\wedge}} \to \pi^0_s(BG).
\]
This statement relates the pure algebraic object of the left hand side to the pure geometric object of the right hand side. The Burnside ring of $G$ is the Grothendieck group of the monoid of isomorphism classes of finite $G$-sets. The sum is induced by disjoint union of $G$-sets and the product by direct product of $G$-sets with diagonal action. The augmentation ideal is the kernel of the augmenation map $A(G)\to \BZ$. Cohomotopy groups are defined as maps to spheres instead of from spheres:
\[
\pi^0_s(BG)=[{\Sigma}^\infty BG_+,\BS]=[\Sigma^\infty BG\vee \BS,\BS].
\]
The map $A(G)\to \pi^0_s(BG)$ (before completion) takes the transitive $G$-set $G/H$ for $H\leq G$ to the map
\[
{\Sigma}^\infty BG_+\stackrel{tr_H}\longrightarrow{\Sigma}^\infty BH_+\to *_+=\BS,
\]
where $tr_H$ is the transfer map. The conjecture was proven by Carlsson, who reduced the case of $p$-groups to the case of $p$-elementary abelian groups. The general case had already been reduced to $p$-groups by McClure. Ravenel, Adams, Guanawardena and Miller used the Adams spectral sequence to do the computations for $p$-elementary abelian groups.

%\subsection{Kervaire invariant problem}
%Look here \url{http://mhovey.web.wesleyan.edu/problems/big.html}

\section{Solutions}
\label{ss:Solutions}

\begin{sol}[Solution to \ref{exe:loopsonS^n}]
\label{exe:loopsonS^nsolution}
This is \cite[p. 243]{Davis-Kirk}.
\end{sol}

\begin{sol}[Solution to \ref{exe:C3C9C3coho}]
\label{exe:C3C9C3cohosolution}
The Lyndon-Hochschild-Serre spectral sequence of the central extension $C_3\to C_9\to C_3$ is 
\[
E_2^{p,q}=H^p(C_3;\BF_3)\otimes H^q(C_3,\BF_3)=\Lambda(y)\otimes \BF_3[x]\otimes \Lambda(y')\otimes \BF_3[x']\Rightarrow H^{p+q}(C_9;\BF_3),
\] 
where $|y|=(1,0),|x|=(2,0),|y'|=(0,1)$ and $|x'|=(0,2)$. So the corner of $E_2$ has the following generators:
\[
\xymatrix@=0pt{
{\phantom{4}} & x'^2 &  x'^2y &  x'^2x & x'^2yx & x'^2x^2\\
{\phantom{3}} & y'x' &  y'x'y &  y'x'x & y'x'yx & y'x'x^2\\
{\phantom{2}} & x'   &  x'y &  x'x & x'yx & x'x^2\\
{\phantom{1}} & y'   &  y'y &  y'x & y'yx & y'x^2\\
{\phantom{0}} & 1    &  y   &  x   & yx   & x^2 \\
{\phantom{0}}\ar@{-}[]+U;[rrrrrr]+UR\ar@{-}[]+R;[uuuuu]+UR  &  {\phantom{0}}  & {\phantom{1}} & {\phantom{2}}   & {\phantom{3}}&{\phantom{4}} &{\phantom{5}}
}
\]
By elementary group theory, $H^1(C_9;\BF_3)=\BF_3$ and hence $y'$ must die killing $d_2(y)=x$ (or $-x$). Then all odd rows disappear in $E_2$ as, for instance,
\begin{align*}
d_2(y'y)=d_2(y')y-y'd_2(y)=xy=yx\text{, }\\
d_2(y'x)=d_2(y')x-y'd_2(x)=x^2.
\end{align*}
Also, we must have $d_2(x')=0$ because $d_2(x')\in \langle y'x\rangle$, $d_2(y'x)=x^2$ and $d_2\circ d_2=0$. Then we deduce that, for instance, 
\[
d_2(x'y)=d_2(x')y'+x'd_2(y)=0,
\]
and in fact $d_2$ is zero on even rows. Summing up,
\[
\xymatrix@=0pt{
{\phantom{4}} & x'^2 &  x'^2y &  x'^2x & x'^2yx & x'^2x^2\\
{\phantom{3}} & y'x'\ar[]+C;[rrd]+C &  y'x'y \ar[]+C;[rrd]+C&  y'x'x\ar[]+C;[rrd]+C & y'x'yx & y'x'x^2\\
{\phantom{2}} & x'   &  x'y &  x'x & x'yx & x'x^2\\
{\phantom{1}} & y'\ar[]+C;[rrd]+C   &  y'y\ar[]+C;[rrd]+C &  y'x \ar[]+C;[rrd]+C& y'yx & y'x^2\\
{\phantom{0}} & 1    &  y   &  x   & yx   & x^2 \\
{\phantom{0}}\ar@{-}[]+U;[rrrrrr]+UR\ar@{-}[]+R;[uuuuu]+UR  &  {\phantom{0}}  & {\phantom{1}} & {\phantom{2}}   & {\phantom{3}}&{\phantom{4}} &{\phantom{5}}
}
\]
$E_3$ is as follows,
\[
\xymatrix@=0pt{
{\phantom{4}} & x'^2 &  x'^2y & 0 & 0 & 0\\
{\phantom{4}} & 0   &  0&  0 & 0 & 0\\
{\phantom{2}} & x'   &  x'y &  0 & 0 & 0\\
{\phantom{1}} & 0   &  0&  0 & 0 & 0\\
{\phantom{0}} & 1    &  y   &  0   & 0   & 0 \\
{\phantom{0}}\ar@{-}[]+U;[rrrrrr]+UR\ar@{-}[]+R;[uuuuu]+UR  &  {\phantom{0}}  & {\phantom{1}} & {\phantom{2}}   & {\phantom{3}}&{\phantom{4}} &{\phantom{5}}
}
\]
and the spectral sequence collapses at the free graded commutative algebra $E_3=E_\infty=\Lambda(y)\otimes \BF_3[x']$. Then by Theorem \ref{thm:freebigradedfreegraded}, we have $H^*(C_9;\BF_3)=\Lambda(z)\otimes \BF_3[z']$ with $|z|=1$ and $|z'|=2$.
\end{sol}

\begin{sol}[Solution to \ref{exe:cohoK(Z,2)}]
\label{exe:cohoK(Z,2)solution}
This is \cite[Example 1.15]{SSch1} and the arguments are already described in \ref{example:cohoK(Z,n)}.
\end{sol}

\begin{sol}[Solution to \ref{exe:cohoK(Z,3)}]
\label{exe:cohoK(Z,3)solution}
This is \cite[Example 1.19]{SSch1} or \cite[p. 245]{Bott-Tu}.
\end{sol}

\begin{sol}[Solution to \ref{exe:coholoopsS^k}]
\label{exe:coholoopsS^ksolution}
This is \cite[Example 1.16]{SSch1} or \cite[p. 204]{Bott-Tu}
\end{sol}

\begin{sol}[Solution to \ref{exe:Poinacreseriesofdihedral}]
\label{exe:Poinacreseriesofdihedralsolution}
The answer is 
\[
P(t)=1+2t+3t^2+4t^3+\ldots=\frac{1}{(t-1)^2}.
\]
\end{sol}

\begin{sol}[Solution to \ref{exe:dihedralasC4:C2}]
\label{exe:dihedralasC4:C2solution}
The Lyndon-Hochschild-Serre spectral sequence is:
\[
H^p(C_2;H^q(C_4;\BF_2))\Rightarrow H^{p+q}(D_8;\BF_2).
\]
The cohomology ring $H^*(C_4;\BF_2)=\Lambda(y)\otimes \BF_2[x]$ with $|y|=1$, $|x|=2$ was described in Example \ref{example:C2C4C2F2}. On each degree $q$, the $\BF_2$-module $H^q(C_4;\BF_2)$ is equal to $\BF_2$, and hence it must be trivial as a $C_2$-module. So, although the given extension is not central, the $E_2$-page is still equal to $E_2^{*,*}=H^*(C_4;\BF_2)\otimes H^*(C_2;\BF_2)=\Lambda(y)\otimes \BF_2[x]\otimes \BF_2[z]$, where $z$ is a generator of $H^1(C_2;\BF_2)$. Generators in the corner of $E_2$ lie as follows:
\[
\xymatrix@=0pt{
{\phantom{4}} & x^2 & x^2z&  x^2z^2& x^2z^3\\
{\phantom{3}} & yx  & yxz &  yxz^2 & yxz^3\\
{\phantom{2}} & x   & xz  &  xz^2  & xz^3\\
{\phantom{1}} & y   & yz  &  yz^2  & yz^3\\
{\phantom{0}} & 1   &  z  &  z^2   & z^3 \\
{\phantom{0}}\ar@{-}[]+U;[rrrr]+UR\ar@{-}[]+R;[uuuuu]+UR  &  {\phantom{0}}  & {\phantom{1}} & {\phantom{2}}   & {\phantom{3}}&{\phantom{4}} &{\phantom{5}}
}
\]
What are the differentials $d_2, d_3, \ldots$? It is straightforward that the Poincar\'e series of $E_2^{*,*}$ is 
\[
P(t)=1+2t+3t^2+4t^3+\ldots=\frac{1}{(t-1)^2}.
\]
This coincides with the Poincar\'e series of $H^*(D_8;\BF_2)$ by Solution \ref{exe:Poinacreseriesofdihedralsolution}. So all terms in the $E_2$-page must survive, all differentials must be zero and $E_\infty=E_2=\Lambda(y)\otimes \BF_2[x]\otimes \BF_2[z]$. By Theorem \ref{thm:Carlsonliftinggensandrels}, $H^*(D_8;\BF_2)=\BF_2[z,\tau,\tau']/(R)$ where $z,\tau,\tau'$ are lifts of $z,y,x$ respectively (see Remark \ref{rmk:horizontalaxisissubalgebraLHS}) and $R$ is a lift of the relation $y^2=0$. This lift must be of the form 
\[
\tau^2=\lambda \tau z+\lambda' z^2,
\]
and we know that $\lambda=1$ and $\lambda'=0$.
\end{sol}

\begin{sol}[Solution to \ref{exe:C2C4C2F2finishoff}]
\label{exe:C2C4C2F2finishoffsolution}
Recall that we have  $E_{\infty}=\BF_2[x,x']/(x^2)=\Lambda(x)\otimes \BF_2[x']$ with $|x|=(1,0)$ and $|x'|=(0,2)$. Hence, by Theorem \ref{thm:Carlsonliftinggensandrels}, we have 
$H^*(C_4;\BF_2)=\BF_2(z,z')/(R)$ with $|z|=1$, $|z'|=2$ and where $R$ is a lift of the relation $x^2=0$. Note that $x^2\in E_\infty^{2,0}=0=F^2H^2\subset H^2(D_8;\BF_2)$ by \ref{equ:filtrationcohoonetotaldimension}. If $z$ is a lift of $x$, then the lift of the relation $x^2=0$ must be $z^2=0$ and we are done.
\end{sol}

\end{document}